\documentclass[12pt,a4paper]{amsart}
\usepackage{amsmath, amssymb, amsthm,
amsfonts, amsthm, bm, eucal, graphicx}
\usepackage[all]{xy}
\usepackage{graphics}
\usepackage{epsfig}
\usepackage{graphicx}
\theoremstyle{plain}
\usepackage{amssymb}
\usepackage[colorlinks]{hyperref}
\usepackage{enumerate}

\advance\hoffset-20mm \advance\textwidth40mm

\newtheorem{thm}{Theorem}[section]
\newtheorem{cor}[thm]{Corollary}
\newtheorem{lem}[thm]{Lemma}
\newtheorem{prop}[thm]{Proposition}

\theoremstyle{definition}
\newtheorem{defi}[thm]{Definition}
\newtheorem{defis}[thm]{Definitions}
\newtheorem{conj}[thm]{Problem}
\newtheorem{conv}[thm]{Convention}
\newtheorem{nota}[thm]{Notation}
\newtheorem{rem}[thm]{Remark}
\newtheorem{rems}[thm]{Remarks}
\newtheorem{exa}[thm]{Example}
\newtheorem{exas}[thm]{Examples}
\newtheorem{sit}[thm]{}

\newcommand{\brem}{\begin{rem}}
\newcommand{\brems}{\begin{rems}}
\newcommand{\erem}{\end{rem}}
\newcommand{\erems}{\end{rems}}
\newcommand{\bexa}{\begin{exa}}
\newcommand{\bexas}{\begin{exas}}
\newcommand{\eexa}{\end{exa}}
\newcommand{\eexas}{\end{exas}}
\newcommand{\bdefi}{\begin{defi}}
\newcommand{\edefi}{\end{defi}}
\newcommand{\bdefis}{\begin{defis}}
\newcommand{\edefis}{\end{defis}}
\newcommand{\bcor}{\begin{cor}}
\newcommand{\ecor}{\end{cor}}
\newcommand{\blem}{\begin{lem}}
\newcommand{\elem}{\end{lem}}
\newcommand{\bconv}{\begin{conv}}
\newcommand{\econv}{\end{conv}}
\newcommand{\bconj}{\begin{conj}}
\newcommand{\econj}{\end{conj}}
\newcommand{\bprop}{\begin{prop}}
\newcommand{\eprop}{\end{prop}}
\newcommand{\bthm}{\begin{thm}}
\newcommand{\ethm}{\end{thm}}
\newcommand{\bnota}{\begin{nota}}
\newcommand{\enota}{\end{nota}}
\newcommand{\bsit}{\begin{sit}}
\newcommand{\esit}{\end{sit}}
\newcommand{\be}{\begin{equation}}
\newcommand{\ee}{\end{equation}}
\newcommand{\bproof}{\begin{proof}}
\newcommand{\eproof}{\end{proof}}
\def\ba{\begin{array}}
\def\ea{\end{array}}
\def\bea{\begin{eqnarray}}
\def\eea{\end{eqnarray}}

\def\bnum{\begin{enumerate}}
\def\enum{\end{enumerate}}
\newcommand{\la}{\label}

\newtheorem*{theo*}{Theorem}

\theoremstyle{definition}

\newtheorem*{definition*}{Definition}

\def\cC{{\mathcal C}}

\def\cN{{\mathcal N}}

\def\cS{{\mathcal S}}

\def\AA{{\mathbb A}}
\def\AT{{\mathbb A}^1_*}

\def\NN{{\mathbb N}}
\def\ZZ{{\mathbb Z}}
\def\QQ{{\mathbb Q}}
\def\TT{{\mathbb T}}
\def\kk{{\Bbbk}}
\def\CC{{\mathbb C}}
\def\G{{\mathbb G}}
\def\CT{{\mathbb C}^{\times}}

\def\PP{{\mathbb P}}

\def\Oc{{\mathcal{C}}}
\def\On{{\mathcal{N}}}

\def\Aut{\mathop{\rm Aut}}

\def\GL{\mathop{\rm GL}}
\def\SL{\mathop{\rm SL}}
\def\Jonq{\mathop{\rm Jonq}}
\def\JONQ{\mathop{\rm JONQ}}

\def\id{\mathop{\rm id}}

\def\Stab{\mathop{\rm Stab}}

\def\deg{\mathop{\rm deg}}

\def\ord{\mathop{\rm ord}}

\def\Aff{\mathop{\rm Aff}}

\def\ll1{l_{\lambda}^{-1}(1)}
\def\lm1{l_{\mu}^{-1}(1)}

\def\ba{\begin{array}}
\def\ea{\end{array}}
\def\bea{\begin{eqnarray}}
\def\eea{\end{eqnarray}}

\begin{document}
\sloppy
\title[Acyclic curves]
{Acyclic curves and group actions\\ on affine toric surfaces}

\author%
{Ivan\ Arzhantsev and Mikhail\ Zaidenberg}
\address{Department of Algebra, Faculty of Mechanics and Mathematics,
Lomonosov Moscow State University, Leninskie Gory 1, Moscow,
119991, Russia } \email{arjantse@mccme.ru}
\address{Universit\'e Grenoble I, Institut Fourier, UMR 5582
CNRS-UJF, BP 74, 38402 St. Martin d'H\`eres c\'edex, France}
\email{Mikhail.Zaidenberg@ujf-grenoble.fr}

\thanks{This work was done
during a stay of the first author at the Institut Fourier,
Grenoble and of the second author at the Max Planck Institute of
Mathematics, Bonn. We thank both institutions for generous support
and hospitality. The work of the first author was supported by a
Simons Award.}
\begin{abstract} We show that every
irreducible, simply connected curve on a  toric affine surface $X$
over $\CC$ is an orbit closure of a $\G_m$-action on $X$. It
follows that up to the action of the automorphism group $\Aut(X)$
there are only finitely many non-equivalent embeddings of the
affine line $\AA^1$ in $X$. A similar description is given for
simply connected curves in the quotients of the affine plane by
small finite linear groups.  We provide also an analog of the
Jung-van der Kulk theorem for affine toric surfaces, and apply
this to study actions of algebraic groups on such surfaces.
\end{abstract}
\dedicatory{To Masayoshi Miyanishi on occasion of his $70$th
birthday}

\maketitle
\bigskip

{\footnotesize \tableofcontents}


\section*{Introduction}

The geometry of  affine toric surfaces still attracts the
researches\footnote{See, for instance, the recent paper
\cite{GZLMH} on the Hilbert scheme of zero-cycles on such a
surface.}. Every affine toric surface over $\CC$ except for $\AT
\times \AT$ and $\AA^1 \times \AT$, where
$\AT=\AA^1\setminus\{0\}$, is the quotient of $\AA^2$ by a small
finite cyclic subgroup $G\subseteq\GL(2,\CC)$.\footnote{A finite
linear group is {\em small} if it does not contain
pseudoreflections.}

Throughout the paper, by `acyclic curve' we mean a connected and
simply connected complex affine algebraic curve. A classification
of acyclic curves on the affine plane, both irreducible and
reducible, up to the action of the automorphism group of the plane
is well known, see e.g. \cite{AM}, \cite{Mi}, \cite{Su},
\cite{Za-85}, \cite{ZL}; we recall it in subsection \ref{ss-AMS}
below. In section \ref{sec-acyc-curves} we classify acyclic curves
on affine toric surfaces. Similarly as in \cite{Za}, actions of
one-parameter groups play a crucial role in this classification.

Let $\pi \colon \AA^2 \to X=\AA^2/\, G$ be the quotient morphism,
and let $\Oc$ be an irreducible acyclic curve on $X$. Then
$\pi^{*}(\Oc)$ is an acyclic (reducible and non-reduced, in
general) curve on $\AA^2$. We show (see theorems \ref{prop3} and
\ref{prop4}) that applying an appropriate automorphism of the
affine plane we can transform the curve $\pi^{*}(\Oc)$ and the
$G$-action on $\AA^2$ to canonical forms simultaneously.

In subsection~\ref{sec1}, given an acyclic plane curve $C$, we
describe the stabilizer subgroup $\Stab(C)\subseteq\Aut(\AA^2)$ of
all automorphisms which preserve $C$. We use this description in
subsections~\ref{sec2} and \ref{sec3} in order to obtain the
canonical forms of irreducible acyclic curves on an affine toric
surface $X$. We treat separately the cases of  the curves passing
or do not passing through the singular point of $X$. This leads in
subsection~\ref{sec4} to the conclusion that every irreducible
acyclic curve on $X$ is the closure of a non-closed orbit of a
$\G_m$-action on $X$. Furthermore, if such a curve is contained in
the smooth locus $X_{\rm reg}$ then it is smooth and is as well an
orbit of a $\G_a$-action on $X$, hence is included in a
one-parameter family of such curves. We show that any affine toric
surface $X$ possesses only finitely many equivalence classes of
embedded affine lines (see theorem \ref{co2}). This is an analog
of the celebrated Abhyankar-Moh-Suzuki Embedding Theorem, which
says that there is just one class of embeddings of the affine line
in the affine plane (see theorem \ref{thm-AMS} below).

The above description enables us to classify in subsection
\ref{reducible} all reduced simply connected curves on an affine
toric surface, whenever they are irreducible or not.

Section \ref{sec-gr-aut} is devoted to the automorphism groups of
affine toric surfaces. In theorem \ref{thm-JvdK-toric-surf} we
obtain  an analog of the classical Jung-van der Kulk theorem on a
free amalgamated product structure on the group $\Aut(\AA^2)$ (see
theorem \ref{thm-JvdK}). Using this theorem we describe in
theorems \ref{thm-KW-toric} and \ref{thm-rdgrac} (reductive)
algebraic groups acting effectively on affine toric surfaces.

In the final section~\ref{sec-non-toric} we deal with acyclic
curves on a quotient $X=\AA^2/G$ of the affine plane $\AA^2$ by a
nonabelian small finite group $G$. It turns out that the only
irreducible acyclic curves on $X$ are the images of the affine
lines in $\AA^2$ passing through the origin (see theorem
\ref{prop6}). In particular, every such curve is the closure of a
$\G_m$-orbit and passes through the singular point. Since the
family of these curves is preserved under automorphisms, the
automorphism group $\Aut(X)$ is rather poor. Namely, it coincides
with $N(G)/G$, where $N(G)$ stands for the normalizer of $G$ in
$\GL(2,\CC)$ (see theorem \ref{thm-non-abel-aut}). Consequently,
$N(G)$ coincides with the normalizer of $G$ in the full
automorphism group $\Aut(\AA^2)$. As an example, we describe
explicitly the affine lines on the quaternion surface
$X=\AA^2/Q_8$.

When the paper was finished, S.~Kaliman kindly informed us that he
also came, for different purposes, to similar conclusions, but
never wrote them down. The authors thank S.~Kaliman for this
information and interesting discussions.


\section{Preliminaries}\la{sec0}

In this section we gather some well known facts on the geometry of
the affine plane $\AA^2$ over $\CC$ that we need in the sequel.


\subsection{Simply connected plane affine  curves}\la{ss-AMS}
By a {\em curve} we mean  (for short) a complex affine algebraic
curve. A curve $C$ is called {\em acyclic} if
$\pi_0(C)=\pi_1(C)=1$, i.e. $C$ is connected and simply connected.
Two plane curves $C$ and $C'$ will be called {\em equivalent}  if
$C'=\gamma(C)$ for some $\gamma\in\Aut(\AA^2)$. The following
theorems provide canonical forms of acyclic and, more generally,
of simply connected plane curves.

\bthm\la{thm-AMS} {\rm (Abhankar-Moh \cite{AM}, Suzuki \cite{Su})}
Any reduced, irreducible, smooth, acyclic plane curve is
equivalent to the affine line $C_y=\{y=0\}$. Furthermore, if $C$
is parameterized via a map
$\varphi:\AA^1\stackrel{\simeq}{\longrightarrow} C$,
$\varphi(t)=(p(t),q(t))$, where $p,q\in\CC[t]$, then either $\deg
p\, | \deg q$ or $\deg q \,| \deg p$. \ethm

Using Suzuki's formula for the Euler characteristics of the fibers
in a fibration on a smooth affine surface \cite{Su} (see also
e.g.\ \cite{Gur}) it is not difficult to deduce the following
corollary.

\bcor\la{cor-AMS}  Any disconnected, simply connected, reduced
plane curve is equivalent to a union of $r\ge 2$ parallel lines.
\ecor

\bthm\la{thm-LZ} {\rm (Lin-Zaidenberg \cite{ZL})} \bnum\item[(a)]
Any reduced, irreducible, singular, acyclic plane curve is
equivalent to one and only one of the curves
$C_{a,b}=\{y^a-x^b=0\}$, where $1< a< b$ and $\gcd(a,b)=1$.
\item[(b)] Any reduced, simply connected plane curve
is equivalent to a curve given by one of the equations
\be\la{eq-par1} y^{\varepsilon_y}p(x)=0\ee or \be\la{eq-par2}
x^{\varepsilon_x}y^{\varepsilon_y}
\prod_{i=1}^r(y^a-\kappa_ix^b)=0\,,\ee where $\varepsilon_x,
\varepsilon_y\in\{0,1\},\,\,\, p\in\CC[t]$ is a polynomial with
simple roots, $a,b\ge 1$ and $\gcd(a,b)=1$, $r> 0$, and
$\kappa_i\in\CC^\times,\,i=1,\ldots,r,$ are pairwise distinct.
\enum\ethm

\subsection{The automorphism group of the affine plane}\la{ssec-Aut}
In this subsection the base field $\kk$ can be  arbitrary. We let
$\AA^n_\kk$ denote the affine $n$-space over $\kk$ and
$\Aff(\AA^n_\kk)$ the group of all affine transformations of this
space. By $\JONQ^+(\AA^2_\kk)$ ($\JONQ^-(\AA^2_\kk)$,
respectively) we denote the group of the de Jonqi\`eres
transformations \be\label{J-plus}\Phi^+:(x,y)\longmapsto (\alpha
x+f(y),\beta y+\gamma) \,,\ee respectively,
\be\label{J-minus}\Phi^-:(x,y)\longmapsto (\alpha x+\gamma,\beta
y+f(x)) \,,\ee where $\alpha,\beta\in\kk^\times$, $\gamma\in\kk$,
and $f\in\kk[t]$. The subgroup \be\la{aff-pm}
{\Aff}^{\pm}(\AA^2_\kk)=
\Aff(\AA^2_\kk)\cap{\JONQ}^{\pm}(\AA^2_\kk)\ee consists of all
upper (lower, respectively) triangular affine transformations
$\Phi^\pm$ with $\deg f\le 1$.
 The structure of the automorphism group $\Aut(\AA^2_\kk)$
 is described by the following  classical
theorem.

\bthm\la{thm-JvdK} {\rm (Jung \cite{Ju}, van der Kulk \cite{vdK})}
The automorphism group $\Aut(\AA^2_\kk)$ is the free product of
the subgroups $\JONQ^+(\AA^2_\kk)$ and $\Aff(\AA^2_\kk)$
amalgamated over their intersection ${\Aff}^+(\AA^2_\kk)$:
$$\Aut(\AA^2_\kk)={\JONQ}^+(\AA^2_\kk)*_{{\Aff}^+(\AA^2_\kk)}
\Aff(\AA^2_\kk)\,.$$ \ethm

\brems\la{hist-rem} 1. In fact Jung \cite{Ju} just established,
over a field $\kk$ of characteristic $0$, the equality
$$\Aut(\AA^2_\kk)=\langle U^+,\,
\Aff(\AA^2_\kk)\rangle\,,$$ where \be\la{max-to} U^+=
\{\Phi\in{\JONQ}^{+}(\AA^2_\kk)\,|\,\Phi:(x,y)\mapsto
(x+f(y),y)\}\,.\ee However, $\Aut(\AA^2_\kk)\neq U^+ *_{
U^+\cap\Aff(\AA^2_\kk)} \Aff(\AA^2_\kk)$, see Remark in \cite[\S
2]{Ka2}.

Over an arbitrary ground field, van der Kulk did not formulate the
theorem in terms of amalgamated free products, but from his
results the theorem can be deduced readily, as this is done in
\cite{Na} or \cite[Theorem 2]{Ka2}. In characteristic zero, the
theorem is a consequence of the Abhyankar-Moh-Suzuki theorem
\ref{thm-AMS}. See also \cite{Al}, \cite{CK}, \cite{Di},
\cite{FB}, \cite{DG}, \cite{Gut1}, \cite{La}, \cite{ML},
\cite{MW}, \cite{NVC}, \cite{Re}, \cite{VdE}, and \cite{Wr3} for
different approaches.

2. Actually we have
$$\Aut(\AA^2_\kk)=\langle{\JONQ}^+(\AA^2_\kk),\,\tau\rangle\,,$$
where $\tau\in\Aff(\AA^2_\kk)$, $\tau:(x,y)\longmapsto (y,x)$ is a
twist. Notice that $\tau\in\langle{\Aff}^+(\AA^2_\kk),\,
{\Aff}^-(\AA^2_\kk)\rangle$ and
$${\JONQ}^-(\AA^2_\kk)
=\tau{\JONQ}^+(\AA^2_\kk)\tau,\quad {\Aff}^-(\AA^2_\kk)
=\tau{\Aff}^+(\AA^2_\kk)\tau,\quad\text{and}\quad U^-=\tau
U^+\tau\,,$$ where the subgroup $U^-\subseteq
{\JONQ}^-(\AA^2_\kk)$ is defined similarly as $U^+$. In particular
$$\Aut(\AA^2_\kk)=\langle{\JONQ}^+(\AA^2_\kk),\,
{\JONQ}^-(\AA^2_\kk)\rangle\,.$$ \erems

The following theorem absorbed several previously known results.
In this generality, it was first proved by Kambayashi \cite{Ka1}
(using a result of Wright \cite{Wr1}) as a consequence of the
Jung-van der Kulk theorem \ref{thm-JvdK}.

\bthm\la{thm-KW} {\rm (Kambayashi \cite[Theorem 4.3]{Ka1}, Wright
\cite{Wr1}, \cite{Wr2})} Any algebraic subgroup of the group
$\Aut(\AA^2_\kk)$ is conjugate either to a subgroup of
$\Aff(\AA^2_\kk)$, or to a subgroup of $\JONQ^+(\AA^2_\kk)$. \ethm

The proof exploites the following observation: every algebraic
subgroup of $\Aut(\AA^2_\kk)$ has bounded degree, hence also a
bounded length with respect to the  free amalgamated product
structure. However, by Serre \cite{Se} a subgroup of bounded
length in an amalgamated free product $A*_C B$ is conjugate to a
subgroup of one of the factors $A$ and $B$.

In the next corollary we suppose  that the base field $\kk$ is
algebraically closed of characteristic $0$ (while certain
assertions remain valid in the positive characteristic case).

\bcor\la{cor-KW} \bnum\item {\rm (Igarashi \cite{Ig}, Furushima
\cite{Fur})} Every finite subgroup of $\Aut(\AA^2_\kk)$ is
conjugate to a subgroup of $\GL(2,\kk)$.

\item  {\rm (Gutwirth \cite{Gut2}, Bialynicki-Birula \cite{BB})}
Every maximal torus in $\Aut(\AA^2_\kk)$ has rank 2 and is
conjugate to the standard maximal torus $\TT\subseteq \GL(2,\kk)$.

\item {\rm (Gutwirth \cite{Gut2}, Bialynicki-Birula \cite{BB})}
Every one-torus in $\Aut(\AA^2_\kk)$  is conjugate to a subtorus
of $\TT$.

\item {\rm (Rentschler \cite{Re})} Every $\G_a(\kk)$-action on $\AA^2_\kk$
is conjugate to an action via de Jonqi\`eres transformations
$$t.(x,y)=(x+tf(y),y),\quad\text{where}\quad t\in \G_a(\kk)
\quad\text{and}\quad f\in\kk[y]\,.$$ \enum\ecor

\brems\la{rem-fixed} 1. Assertions (1)-(3) follow from a more
general result for reductive groups, see \ref{prop-quasitorus}
below.

2. Analogs of (1) and (4) fail in higher dimensions, while (3)
holds in dimension $3$ and is open in higher dimensions. We do not
dwell on this here (see, however, \cite{KMLKR}, \cite{Po} and the
survey \cite{Kr}; see also \cite{As} and \cite{MN} for the case of
a positive characteristic). \erems


\section{Subgroups of de Jonqu\`eres group
and  stabilizers of plane curves}


\subsection{Subgroups of the de Jonqu\`eres group} In this subsection
the base field $\kk$ will be an algebraically closed field of
characteristic zero (while some results are still valid over an
arbitrary field of characteristic zero.) By abuse of language, we
still call {\em de Jonqu\`eres groups} the subgroups
$\Jonq^\pm(\AA^2_\kk)\subseteq \JONQ^\pm(\AA^2_\kk)$ consisting,
respectively, of the transformations \be\la{eq-Jonq-plus}
\varphi^+:(x,y)\longmapsto (\alpha x+f(y),\beta y) \,,\ee and
\be\la{eq-Jonq-minus} \varphi^-:(x,y)\longmapsto (\alpha x,\beta
y+f(x)) \,,\ee where $\alpha,\beta\in\kk^\times$ and $f\in\kk[t]$
(so $\Phi^\pm$ as in (\ref{J-plus}), (\ref{J-minus}) belongs to
$\Jonq^\pm(\AA^2_\kk)$ if and only if $q=0$).

Clearly, \be\la{eq-unip} {\Jonq}^{+}(\AA^2_\kk)\simeq
U^+\rtimes\TT\,,\ee where $\TT\subseteq\GL(2,\kk)$ is the maximal
torus
$$ \TT=\{\delta\in{\Jonq}^+(\AA^2_\kk)\,|\,\delta:(x,y)\mapsto
(\alpha x, \beta y)\}\,,$$ and $U^+\simeq\kk[t]$ as in
(\ref{max-to})
 is an infinite dimensional vector
group. Indeed, let $\rho:{\Jonq}^{+}(\AA^2_\kk)\to\TT$ denote the
canonical surjection provided by (\ref{eq-unip}). Then any element
$\varphi^+\in{\Jonq}^{+}(\AA^2_\kk)$
 as in (\ref{eq-Jonq-plus})
 admits a decomposition $\varphi^+=\partial\circ\delta$,
 where $\delta=\rho(\varphi^+)\in\TT$
 and $\partial\in U^+$,
$\partial:(x,y)\mapsto (x+f(y/\beta),y)$.

In particular, ${\Jonq}^{+}(\AA^2_\kk)$ is a metabelian group, and
$U^+$ can be considered as its unipotent radical in the following
sense. There is a natural filtration by algebraic subgroups
$${\Jonq}^{+}(\AA^2_\kk)
=\bigcup_{n\in\NN}{\Jonq}^{+}_{n}(\AA^2_\kk)\,,$$ where
$${\Jonq}^{+}_{n}(\AA^2_\kk)\simeq U^+_n\rtimes\TT\,$$ consists of
all elements $\varphi^+\in{\Jonq}^{+}(\AA^2_\kk)$ with $\deg f\le
n$. Then $\TT$ is a maximal torus of ${\Jonq}^{+}_{n}(\AA^2_\kk)$,
and $U^+_n=U^+\cap {\Jonq}^{+}_{n}(\AA^2_\kk)$ is its unipotent
radical.

Any algebraic subgroup $G$ of ${\Jonq}^{+}(\AA^2_\kk)$ has finite
degree
$$d= \min\{n\,|\,G\subseteq {\Jonq}^{+}_n(\AA^2_\kk)\}\,.$$
In particular, every maximal torus in ${\Jonq}^{+}(\AA^2_\kk)$ is
conjugate with $\TT$, and any algebraic subgroup which consists of
semi-simple elements is conjugate to a subgroup of $\TT$ (see
\cite[VII.19.4, \, VIII.21.3A]{Hum}).

In the following lemma we characterize semi-simple and torsion
elements in the group ${\Jonq}^{+}(\AA^2_\kk)$.

\blem\la{lem-uni-rad}\bnum\item[(a)]  An element
$\varphi^+\in{\Jonq}^{+}(\AA^2_\kk)$ as in (\ref{eq-Jonq-plus}) is
semi-simple if and only if the polynomial $f(y)=\sum_{m\ge 0}
a_my^m$ satisfies the condition \be\la{eq-fin-ord} a_m=0
\quad\text{if}\quad \alpha=\beta^m\,,\ee if and only if there
exists $\mu\in U^+$ such that $\mu^{-1}\varphi^+\mu=\delta\in\TT$
\footnote{In fact here $\delta=\rho(\varphi^+)$.}.
\item[(b)] An element $\varphi^+
\in {\Jonq}^{+}(\AA^2_\kk)$ as in (\ref{eq-Jonq-plus}) is of
finite order if and only if it is semi-simple and
$\delta=\rho(\varphi^+)\in\TT$ is of finite order. \enum \elem

\bproof (a) We claim that if $\varphi^+$ satisfies
(\ref{eq-fin-ord}) then there exists $\mu\in U^+$,
$\mu:(x,y)\mapsto (x+g(y),y)$, where $g\in\kk[t]$, such that
\be\la{eq-commuta}
\mu^{-1}\varphi^+\mu=\delta=\rho(\varphi^+)\in\TT\quad\text{or,
equivalently,}\quad\partial=[\mu,\delta]=\mu\delta\mu^{-1}\delta^{-1}\,.
\ee Writing $g(y)=\sum_{m\ge 0} b_my^m$, it is readily seen that
$\mu$ satisfies (\ref{eq-commuta}) if and only if  the
coefficients $b_m$ of $g$ satisfy the conditions \be\la{eq-syst}
b_m=\begin{cases}\frac{a_m}{\beta^m-\alpha} &\text{ if
$\beta^m-\alpha\neq 0$}\\ \text{arbitrary} &\text{if
$\beta^m-\alpha= 0=a_m$.}\end{cases} \ee Indeed,
(\ref{eq-commuta}) can be written as
$$ (x+g(y)-\alpha g(y/\beta),y)=(x+f(y/\beta),y)\,$$
which is equivalent to (\ref{eq-syst}). This shows the existence
of $\mu$ in (\ref{eq-commuta}) under condition (\ref{eq-fin-ord}).

On the other hand, if (\ref{eq-fin-ord}) fails i.e., there exists
$m\in\NN$ such that $\beta^m-\alpha= 0$ and $a_m\neq 0$, then
$\mu$ as in (\ref{eq-commuta}) cannot exist. The remaining claims
in (a) are easy and so we leave them to the reader.

(b) We have to show that if $(\varphi^+)^k=\id$ then
(\ref{eq-fin-ord}) holds and $\delta^k=\id$, and vice versa.
Indeed, letting $\gamma=\varphi^+$ for any $k\ge 1$ we obtain
$$\gamma^k:(x,y)\longmapsto
\left(\alpha^kx+f_k(y),\, \beta^k y\right)\,,$$ where
$$f_k(y)=\alpha^{k-1}f(y)+\alpha^{k-2}f(\beta y)+\ldots+
f(\beta^{k-1} y)\,.$$ Therefore $\gamma^k=\id$ if and only if
$\alpha^k=\beta^k=1$ (i.e.\ $\delta^k=\id$) and $f_k=0$. However,
$f_k=0$ if and only if
$$\forall m\ge 0,\qquad \text{either}\quad a_m=0\quad\text{or}\quad
\alpha^{k-1} + \alpha^{k-2}\beta^m + \ldots + \beta^{m(k-1)}=0\,.
$$
The latter equality can be written as
$$
\left(\alpha/\beta^m\right)^{k-1} +
\left(\alpha/\beta^m\right)^{k-2} + \ldots + 1=0\,.
$$
Since $\left(\alpha/\beta^m\right)^k=1$ this holds if and only if
$\alpha\neq\beta^m$. \eproof

Lemma \ref{lem-uni-rad}(b) admits the following  interpretation.

\brem\la{rem-grading} Consider a $\ZZ/d\ZZ$-grading on the
polynomial ring $A=\kk[t]$:
$$A=\bigoplus_{i=0}^{d-1} A_{d,i},\quad\text{where}\quad
A_{d,i}=t^i\kk[t^d]\,.$$ According to this decomposition any
element $f\in A$ can be written as $f=\sum_{i=0}^{d-1} f_i$, where
$f_i\in A_i$ $\forall i$. Assuming that $\beta$ has finite order
$d$ we let $m(\alpha,\beta)=\min\{m\ge 0\,|\,\alpha=\beta^m\}$.
Then condition (\ref{eq-fin-ord}) can be expressed as
$f_{m(\alpha,\beta)}=0$, provided this quantity is well defined.
(Otherwise (\ref{eq-fin-ord}) does not impose any condition.) This
phenomenon can be seen on the following simple examples. \erem

\bexa\la{ex-finite-order} Letting $d=2$ any element $f\in\kk[t]$
can be written as $f=f_0+f_1$, where $f_0$ is even and $f_1$ is
odd. There are the following three types of involutions
$\varphi^+:(x,y)\mapsto (\alpha x+f(y),\beta y)$ in
${\Jonq}^{+}(\AA^2_\kk)$: \bnum\item $\alpha=1,\,\,\,\beta=-1$,
$f\in\kk[t]$ is odd i.e. $f_0=0$;
\item $\alpha=-1,\,\,\,\beta=1$, $f\in\kk[t]$ is arbitrary;
\item $\alpha=-1,\,\,\,\beta=-1$, $f\in\kk[t]$ is even i.e. $f_1=0$.\enum
\eexa

\blem\la{lem-commuting} Consider a pair of elements
$\gamma,\,\tilde{\gamma}\in{\Jonq}^{+}(\AA^2_\kk)$,
$$\gamma:(x,y)\mapsto (\alpha x+f(y),\beta y)\quad\text{and}\quad
\tilde{\gamma}:(x,y)\mapsto (\tilde{\alpha}
x+\tilde{f}(y),\tilde{\beta} y)\,,$$ where
$$f(y)=\sum_{m\ge 0}a_my^m\quad\text{and}\quad \tilde{f}(y)=\sum_{m\ge
0}\tilde{a}_my^m\,.$$ Then $\gamma$ and $\tilde{\gamma}$ commute
if and only if \be\la{eq-commuting}
a_m({\tilde{\beta}}^{m}-\tilde{\alpha})=\tilde{a}_m(\beta^{m}-\alpha)\quad\forall
m\ge 0\,.\ee \elem

\bproof The proof is easy and is left to the reader. \eproof

Recall that a {\em quasitorus} is a product of a torus and a
finite abelian group. Any algebraic subgroup of a torus is a
quasitorus.

\bprop\la{prop-quasitorus} Any reductive algebraic subgroup $G$ of
the group ${\Jonq}^{+}(\AA^2_\kk)$ is
conjugate to a subgroup of the torus $\TT$. More precisely, there
exists an element $\mu\in U^+$ such that \be\la{eq-mu}
\mu^{-1}\circ\varphi\circ\mu= \rho(\varphi)\quad\forall\varphi\in
G\,, \ee where $$\rho:{\Jonq}^{+}(\AA^2_\kk)\to\TT=
{\Jonq}^{+}(\AA^2_\kk)/U^+$$ is the natural surjection. In
particular $\mu^{-1}G\mu= \rho(G)\subseteq\TT\,.$ \eprop

\bproof Since $G$ is reductive the unipotent radical $R_u(G)$ is
trivial. Hence also the subgroup $G\cap U^+$ is trivial. Thus
$\rho$ restricts to an isomorphism
$\rho|G:G\stackrel{\simeq}{\longrightarrow} \rho(G)\subseteq\TT$.
In particular, $G$ is abelian and
 consists of semi-simple elements, cf.\ lemma \ref{lem-uni-rad}. By
\cite[VII.19.4, \, VIII.21.3A]{Hum} $G$ is contained in a maximal
torus  $\TT'$ conjugate to $\TT$. Now the first assertion follows.

Let us show the second. Since $G$ is abelian  and
 consists of semi-simple elements, (\ref{eq-fin-ord}) and
(\ref{eq-commuting}) are fulfilled for any pair of elements
$\gamma,\tilde{\gamma}\in G$. Thus there is $\mu\in U^+$
satisfying (\ref{eq-syst}) and then also (\ref{eq-commuta}) for
all $\varphi\in G$ simultaneously (see the proof of lemma
\ref{lem-uni-rad}). This $\mu$ is as desired. \eproof

As an application we can deduce the following well known fact (see
\cite[Theorem 2]{KS}).

\bcor\la{cor-lineariz-redgr} Every effective action of a reductive
algebraic group $G$ on the affine plane $\AA^2_\kk$ is
linearizable. In other words, $G$ is conjugate in the group
$\Aut(\AA^2)$ to a subgroup of  $\GL(2,\kk)$. \ecor

Let us provide an argument following an indication in \cite{KS}.

\bproof By theorem \ref{thm-KW} $G$ is conjugate in  $\Aut(\AA^2)$
to a subgroup of one of the groups $\Aff(\AA^2_\kk)$ and
$\Jonq^+(\AA^2_\kk)$. In the latter case by proposition
\ref{prop-quasitorus} $G$ is conjugate to a subgroup of the torus
$\TT$ contained in $\Aff(\AA^2_\kk)$. Thus we may assume that
$G\subseteq\Aff(\AA^2_\kk)$.

It remains to show that $G$ has a fixed point in $\AA^2$. Observe
that $G\subseteq\Aff(\AA^2_\kk)$ admits a representation in
$\GL(3,\kk)$ by matrices of the form $\begin{pmatrix}
* & * & * \\
* & * & * \\
0 & 0 & 1 \\
\end{pmatrix}$. Since $G$ is reductive and char$\,(\kk)=0$, $G$
is geometrically reductive. Therefore the $G$-invariant plane
$L_0=\{x_3=0\}$ in $\AA^3_\kk$ with coordinates $(x_1,x_2,x_3)$
has a $G$-invariant complement, say, $R$, which meets the parallel
$G$-invariant plane $L_1=\{x_3=1\}$ in a fixed vector. This yields
a fixed point of $G$ in $\AA^2_\kk\simeq_G L_1$. Now the proof is
completed.
 \eproof

\brem\la{rem-linear} Due to corollary \ref{cor-lineariz-redgr},
when dealing with quotients of the affine plane by finite group
actions it suffices to restrict to linear such actions. \erem

\subsection{Stabilizers of  acyclic plane curves}
\label{sec1}

In the sequel $\AA^n$ stands for the affine $n$-space over $\CC$.
For an algebraic curve $C$ in $\AA^2$ we let $\Stab(C)$ be the
stabilizing group, or {\em stabilizer} of $C$ i.e. the group of
all algebraic automorphisms of $\AA^2$ that preserve $C$:
$$
\Stab(C) \ = \ \{ \gamma \in \Aut (\AA^2) \, | \, \gamma(C)=C\}\,.
$$
In this subsection we describe the stabilizers of acyclic plane
curves given in one of the canonical forms (\ref{eq-par1}) and
(\ref{eq-par2}) of theorem \ref{thm-LZ}(b). Accordingly, we
distinguish the following six types of acyclic curves $C$:

\smallskip

\bnum\item[(I)] $C$ is irreducible and smooth, and then $C\sim
C_y$; \item[(II)] $C$ consists of two smooth, mutually transversal
components, and then $C\sim C_x\cup C_y$;
\item[(III)] $C$ has $r\ge 2$ singular points
and so is equivalent to a union of $r$ parallel lines and a
transversal line given by (\ref{eq-par1}) with $\varepsilon_y=1$;
\item[(IV)] $C$ has an ordinary singularity of multiplicity $r\ge
3$  and so is equivalent to a union of $r$ distinct lines through
the origin;
\item[(V)] All irreducible components of $C$ are smooth and $C$
has a non-ordinary singular point. So  $C$ is equivalent to a
curve given by equation (\ref{eq-par2}) with $a=1<b$ and
$\varepsilon_y+r\ge 2$; \item[(VI)] $C$ contains a singular
component,  and then it is equivalent to a curve given by
(\ref{eq-par2}) with $\min\{a,b\}>1$ and $r\ge 1$. \enum We
analyse each of these cases separately.
\smallskip

In the next proposition we study the curves of type (I). We show
that  the stabilizer $\Stab(C_y)$ of the coordinate axis
$C_y=\{y=0\}$ in $\AA^2$ consists of the de Jonqui\`eres
transformations. For the acyclic curves of types (II)-(IV) reduced
to the canonical form the stabilizer is described in corollaries
\ref{cross-lin} and \ref{cross}, and for those of types (V) and
(VI) in propositions \ref{cross-second} and \ref{cross-third}.

\bprop\label{prop1} The stabilizer $\Stab(C_y)$ in $\Aut(\AA^2)$
coincides with the subgroup $\Jonq^+(\AA^2)$, while
$\Stab(C_x)=\Jonq^-(\AA^2)$. \eprop

\bproof Every $\gamma\in\Stab(C_y)$ sends $y$ to $\beta y$ for
some $\beta\in\CT$. Up to an affine transformation we may assume
that $\beta = 1$ and $\gamma |_{C_y} = \id_{C_y}$. Suppose that
$\gamma$ sends $x$ to $x+h(x,y)$. Since $h(x,0)=0$ we have
$h(x,y)=yp(x,y)$. To show that $p(x,y)$ does not depend on $x$ we
write
$$
p(x,y) \, = \, a_0(y)+a_1(y)x+\ldots+a_k(y)x^k \quad \text{with}
\quad a_k(y) \ne 0\,.
$$
Clearly, $\gamma$ preserves every line $y=y_0$ and induces an
affine automorphism of this line. Picking $y_0$ with $a_k(y_0)\ne
0$ we get $k\le 1$. Letting $k=1$ we obtain $\gamma:x\mapsto
(1+ya_1(y))x+ya_0(y)$. If $y_1$ is  a root of the  non-constant
polynomial $1+ya_1(y)$ then $\gamma$ induces a constant map on the
line $y=y_1$, a contradiction. Hence $k=0$. Thus
$\Stab(C_y)=\Jonq^+(\AA^2)$. The proof of the second assertion is
similar. \eproof

In the following two corollaries we describe the stabilizers of
the canonical curves of types (II)--(IV).

\bcor \label{cross-lin} If $C=\bigcup_{i=1}^r L^i$ is a union of
$r\ge 2$ affine lines in $\AA^2$ through the origin then
$\Stab(C)\subseteq\GL(2,\CC)$.\ecor

\bproof We may suppose that $L_1=C_x$ and $L_2=C_y$. For any $g\in
\Stab(C)$ we can find $h\in\GL(2,\CC)$ such that $g(L_i)=h(L_i)$,
$i=1,2$. It follows by proposition \ref{prop1} that
$\gamma=h^{-1}g\in\GL(2,\CC)$. Hence also
$g=h\gamma\in\GL(2,\CC)$. \eproof

The following corollary is immediate.

\bcor \label{cross} \bnum
\item[(a)] Let $C=\{f(y)=0\}$, where $f\in\CC[y]$
is a polynomial of degree $\ge 2$ with simple roots. If
$K\subseteq\CC$ denotes the set of these roots, then
$$\Stab(C)=\TT_{1,0}\cdot U^+\cdot\Stab(K)\,,$$ where
$U^+\subseteq\Aut(\AA^2)$ is as in (\ref{max-to}),
$$\TT_{1,0}=\{\lambda\in\TT\,|\,\lambda:(x,y)\mapsto (\alpha x,
y),\,\alpha\in\CC^{\times}\}\,,$$ and the stabilizer
$\Stab(K)\subseteq\Aut(\AA^1)\hookrightarrow\Aut(\AA^2)$ acts
naturally on the symbol $y$.
\item[(b)] If $C$ of type (II) is the
coordinate cross $\{xy=0\}$ in $\AA^2$ then $\Stab(C)=N(\TT)$ is
the normalizer of the maximal torus $\TT$ in the group
$\GL(2,\CC)$.
\item[(c)] If $C$ of type (IV) is a union of $r$ affine
lines through the origin, where $r\ge 3$, then
$\Stab(C)\subseteq\GL(2,\CC)$ is  a finite extension of the group
$\TT_{1,1}=\CC^\times\cdot \id$ of scalar matrices.
 \item[(d)] If $C$ of type (III) is given by equation $xyf(x)=0$, where $f\in\CC[x]$
is a polynomial of degree $\ge 1$ with simple roots such that
$f(0)\neq 0$, then $\Stab(C)\subseteq\TT$ is a finite extension of
the one-torus $\TT_{0,1}\subseteq\TT$. \enum \ecor

Thus in (b)-(d) the stabilizer $\Stab(C)$ is a linear group, while
the group in (a)  is infinite dimensional.  The group $\Stab(C)$
in (b) is nonabelian. It can occur to be nonabelian also in (c),
as in the following simple example.

\bexa\label{nonab} Given a  nonabelian finite subgroup
$G\subseteq\GL(2,\CC)$ consider the curve $C=\bigcup_{g\in
G}g(C_y)$. Since $G\subseteq\Stab(C)$ the latter group is
nonabelian.\eexa

Consider further  an irreducible acyclic curve $C_{a,\, b}$ given
in $\AA^2$ by equation $y^a-x^b=0$, where $a,b\ge 1$ and
$\gcd(a,b)=1$. In the following proposition we describe the
stabilizer $\Stab(C_{a,\, b})$ for a singular such curve. Consider
a one-parameter subgroup $\TT_{a,b}$ of the torus $\TT$,
\be\label{eq-1}\TT_{a,b}=\{\gamma_{a,b}(t)\,|\,t\in\G_m\}
\subseteq \Stab(C_{a,\, b}),\quad\text{where}\quad
\gamma_{a,b}(t): (x,y) \longmapsto (t^a x,\,t^b y)\,.\ee

\bprop \label{prop2} If $\min\{a,b\}>1$ then $\Stab(C_{a,\,
b})=\TT_{a,b}$. \eprop

\bproof Letting $C=C_{a,\, b}$ and $\Gamma=\Stab(C)$, we consider
the pointwise stabilizer
$$\Gamma_0= \{ \gamma \in \Gamma \, | \, \gamma|_C = {\id}_C
\}\subseteq \Gamma\,.$$

\smallskip

\noindent {\em Claim 1. The group $\Gamma_0$ is torsion free.}

\smallskip

\noindent {\em Proof of claim 1.} Let $\gamma_0 \neq \id$ be an
element of finite order in $\Gamma_0$. The finite cyclic group
$\langle \gamma_0\rangle$ is reductive. Hence  the fixed point
locus $(\AA^2)^{\gamma_0}$ is smooth by the Luna \'Etale Slice
Theorem~\cite{Lu} (see also \cite[\S6]{PV}). Thus an irreducible
component $C=C_{a,b}$ of this locus must be smooth as well.
However, under our assumptions the curve $C_{a,\, b}$ is singular.
\qed

\smallskip

An element $\gamma\in \Gamma$ sends the polynomial $q(x,y) =
y^a-x^b$ to $\lambda_{\gamma}q$, where $\lambda_{\gamma} \in
\G_m$. Letting $\psi(\gamma)=\lambda_{\gamma}$ yields a character
$\psi \colon \Gamma \to \G_m$ of $\Gamma$. The following claim is
immediate.

\smallskip

\noindent {\em Claim 2.} $\Gamma=\TT_{a,b}\cdot\Gamma_1$, where
$\Gamma_1=\ker (\psi)$.

\smallskip

\noindent {\em Claim 3. The kernel $\Gamma_1=\ker (\psi)$ is a
torsion group. Furthermore, there is a positive integer $N$ such
that the orders of all elements in $\Gamma_1$ divide $N$.}

\smallskip

\noindent {\em Proof of claim 3.} The group $\Gamma_1$ acts on
every fiber  $C_{\alpha}=\{q(x,y)=\alpha\}$ of $q$. Since for
$\alpha\ne 0$ the affine curve $C_{\alpha}$ has positive genus,
its automorphism group is finite of order, say, $N$. For any
$\alpha \ne 0 \ne \beta$ the curves $C_{\alpha}$ and $C_{\beta}$
are isomorphic. So  their automorphism groups are isomorphic, too.
Since $\gamma^N|_{C_{\alpha}}=\id_{C_{\alpha}}$ for every
$\gamma\in\Gamma_1$ and $\alpha\ne 0$, we have $\gamma^N=\id$.
\qed

\smallskip

\noindent {\em Claim 4. } $\Gamma_0=\{\id\}$.

\smallskip

\noindent {\em Proof of claim 4.} According to claim 3 we have
$\Gamma_0\subseteq\TT_{a,b}\cdot\Gamma_1$. Writing an element
$\gamma_0\in\Gamma_0$ as $\gamma_0=\gamma_{a,b}(t)\circ\gamma_1$,
from $\gamma_0|_C=\id_C$ we obtain
$\gamma_1^{-1}|_C=\gamma_{a,b}(t)|_C$. Hence
$\id_C=\gamma_1^{-N}|_C=\gamma_{a,b}(t^N)|_C$. It follows that
$t^N=1$. Since $\Gamma_0\cap\Gamma_1=\{\id\}$ the map
$\psi|_{\Gamma_0}:\Gamma_0\to\G_m$ is injective. So
$\psi(\gamma_0)=\psi(\gamma_{a,b}(t))$ has finite order dividing
$N$. Due to claim 1 we can conclude that $\Gamma_0=\{\id\}$.\qed

\smallskip

\noindent {\em Claim 5. } $\Gamma=\TT_{a,b}$.

\smallskip

\noindent {\em Proof of claim 5.} For any $\gamma\in\Gamma$ there
exists $t\in \CC^\times$ such that $\gamma|C=\gamma_{a,b}(t)|C$.
Hence $\gamma\circ\gamma_{a,b}^{-1}(t)\in\Gamma_0=\{\id\}$ and so
$\gamma=\gamma_{a,b}(t)\in\TT_{a,b}$.\qed

This ends the proof. \eproof

Using propositions \ref{prop1} and \ref{prop2} we describe in
\ref{cross-second} and \ref{cross-third} below the structure of
the stabilizer $\Stab(C)$ for  reduced (but possibly reducible)
acyclic plane curves $C$  of the remaining types (VI) and (V),
respectively.

\bprop\label{cross-second} Let  $C$ be an acyclic  curve  of type
(VI) given by equation (\ref{eq-par2}), where $r\ge 1$ and
$\gcd(a,b)=1$. If $\min\{a,b\}>1$ then $\Stab(C)$ is a quasitorus
of rank 1 contained in the maximal torus $\TT$. \eprop

\bproof Let $C^i=\{y^a-\kappa_ix^b=0\}$, $i=1,\ldots,r$, be the
irreducible components of the curve $C$. Clearly
$\TT_{a,b}\subseteq \Stab(C)$. If $r=1$ then by proposition
\ref{prop2} $\Stab(C)=\TT_{a,b}$. Suppose that $r>1$. Consider a
finite abelian group $H=\Stab(C)\cap\TT_{1,0}$. We claim that
$\Stab(C)=H\cdot\TT_{a,b}\subseteq\TT$. Hence this is a quasitorus
of rank one, as stated. Indeed, if
$\delta\in\Stab(C)\setminus\TT_{a,b}$ then $\delta(C^1)=C^i$ for
some $i>1$. If $h\in\TT_{1,0}$ is such that $h(C^i)=C^1$ then
$\gamma=h\circ\delta\in\Stab(C^1)=\TT_{a,b}$. Hence
$h=\gamma\circ\delta^{-1} \in H$ and so
$\delta=h^{-1}\circ\gamma\in H\cdot\TT_{a,b}$. Now the claim
follows. This ends the proof. \eproof

\bprop\label{cross-third} Let $C$ be an acyclic curve $C$ of type
(V) given by equation (\ref{eq-par2}), where $r\ge 1$ and
$\varepsilon_y+r\ge 2$. If $a=1<b$ then $\Stab(C)$ is a quasitorus
of rank 1 conjugated in the group $\Aut(\AA^2)$ to a subgroup of
the torus $\TT$. \eprop

\bproof If $\varepsilon_x=1$ the proof is easy and can be left to
the reader (cf.\ corollary \ref{cross}(b)). Thus we may  restrict
to the case $\varepsilon_x=0$. By our assumptions $C$ is reducible
and all components $C^i$ of $C$ are smooth and mutually tangent at
the origin (we let here $C^0=C_y$ if $\varepsilon_y=1$). Consider
the pencil $\mathcal L=\{C_\mu\}_{\mu\in\PP^1}$, where
$C_\mu=\{y-\mu x^b=0\}$ for $\mu\neq\infty$ and $C_\infty=C_x$.

\smallskip

\noindent {\em Claim 1. The pencil $\mathcal L\setminus\{C_x\}$ is
stable under the action of the group $\Stab(C)$ on $\AA^2$.}

\smallskip

\noindent {\em Proof of claim 1.} The unique singular point $\bar
0\in C$ is fixed under the action of $\Stab(C)$. Furthermore, for
every $g\in\Stab(C)$ and every $\mu\in\CC$, either $C_\mu\subseteq
C$ or $g(C_\mu)\cap C=\{\bar 0\}$. In the latter case, letting
$B_\mu=g(C_\mu)$ it follows that the restriction $(y-\kappa_i
x^b)|_{B_\mu}$ vanishes just at the origin. Hence in an affine
coordinate, say,  $z$ in $B_\mu\simeq\AA^1$ centered at the origin
the latter function is a monomial $\lambda_i z^{\alpha_i}$, where
 $\lambda_i\in\CC^\times$ and $\alpha_i\ge 1$. Therefore
$B_\mu$ is a component of the curve
$$(y-\kappa_i x^b)^{\alpha_j}-\delta(y-\kappa_j
x^b)^{\alpha_i}=0,\quad\text{where}
\quad\lambda_i^{\alpha_j}=\delta\lambda_j^{\alpha_i}\,.$$ Since
$B_\mu$ is smooth one of the exponents $\alpha_i,\alpha_j$ divides
the other. Preservation of the local intersection indices under
$g$ implies that $i(B_\mu,C^i,\bar 0)=b$ $\forall i$. It follows
that actually $\alpha_i=\alpha_j$ $\forall i,j$. Finally $B_\mu$
coincides with a certain member $C_{\mu'}\in\mathcal L$, as
claimed. \qed

\smallskip

Clearly $\nu(g): \AA^1\to\AA^1$, $\mu\mapsto\mu'$, is an affine
transformation leaving invariant the set
$K=\{\kappa_i\}\subseteq\AA^1$ enriched by $\kappa_0=0$ in case
where $\varepsilon_y=1$ in (\ref{eq-Jonq-plus}). Let $\Stab(K)$ be
the stabilizer of $K$ in $\Aut(\AA^1)$. The group $\Stab(K)$ fixes
the isobarycentre of $K$ and so embeds in $\G_m$ onto a finite
subgroup. The natural homomorphism $\varphi: \Stab(C)\to\Stab(K)$
fits in the exact sequence
$$1\to{\Stab}_0(C)\to \Stab(C)\stackrel{\varphi}
{\longmapsto}\Stab(K)\to\ldots\,,$$ where
$\Stab_0(C)=\ker(\varphi)\subseteq \Stab(C)$ consists of the
elements $g\in\Aut(\AA^2)$ leaving invariant every component $C^i$
of $C$.

\smallskip

\noindent {\em Claim 2. $\Stab_0(C)=\TT_{1,b}$.}

\smallskip

\noindent {\em Proof of  claim 2.} If $g\in\Stab_0(C)$ then
$\varphi(g)=\id_K$. Since $|K|\ge 2$ it follows that
$\nu(g)=\id_{\AA^1}$ i.e., $g(C_\mu)=C_\mu$ $\forall \mu\in\AA^1$.
In particular, $g(C_0)=C_0$, where $C_0=C_y$. By proposition
\ref{prop1} we have $g:(x,y)\mapsto (\alpha x+f(y),\beta y)$ for
some $\alpha,\beta\in\CC^{\times}$ and $f\in\CC[y]$. The equality
$g(C^1)=C^1$ implies that $f=0$ and $\beta=\alpha^b$, that is,
$g\in\TT_{1,b}$. Now the claim follows. \qed

\smallskip

Thus $\Stab(C)$ is an extension of the one-torus $\TT_{1,b}$ by a
finite cyclic group. The proof can be completed due to the
following

\smallskip

\noindent {\em Claim 3. $\Stab(C)$ is conjugated in $\Aut(\AA^2)$
to a subgroup of the maximal torus $\TT$. }

\smallskip

\noindent {\em Proof of  claim 3.}  For every $g\in\Stab(C)$ we
have $g(C_x)=C_x$. Letting $z_0$ denote the isobarycentre of $K$
we consider an automorphism $\gamma\in\Aut(\AA^2)$, $(x,y)\mapsto
(x,y-z_0x^b)$. It is easily seen that $\gamma$ preserves the
pencil $\mathcal L$, while $\nu(\gamma): z\mapsto z-z_0$. Since
$\nu(\gamma
g\gamma^{-1})=\nu(\gamma)\nu(g)\nu(\gamma)^{-1}:0\mapsto 0$ we
obtain $\gamma g\gamma^{-1}:C_y\to C_y$ and $C_x\to C_x$. Hence
$\gamma\circ g\circ\gamma^{-1}\in\TT$ (see corollary
\ref{cross}(b)) and so $\gamma \Stab(C)\gamma^{-1}\subseteq\TT$.
\eproof

However,  for $C$ as in proposition \ref{cross-third} the
stabilizer $\Stab(C)$ is not necessarily contained in
$\GL(2,\CC)$, as is seen in the following exemple.

\bexa\label{ex-non-lin} Let $C=\{y(y-x^b)=0\}$, where $b>1$, and
let $g:(x,y)\mapsto (x,x^b-y)$. Then
$g\in\Stab(C)\setminus\GL(2,\CC)$. \eexa

From  corollary \ref{cross} and propositions
\ref{prop2}-\ref{cross-third} we deduce the following result.

\bcor\label{final-cor}  The stabilizer $\Stab(C)$ of an acyclic
plane curve $C$  is abelian unless $C$ is equivalent under the
$\Aut(\AA^2)$-action to a  union of affine lines through the
origin. \ecor


\section{Acyclic curves on affine toric surfaces}
\label{sec-acyc-curves}

In this section we classify acyclic curves on affine toric
surfaces, similarly as this is done for plane acyclic curves in
theorems \ref{thm-AMS} and \ref{thm-LZ}.


\subsection{Acyclic curves in the smooth locus}
\label{sec2}

We start with the curves that do not pass through the singular
point.

Let $d$ and $e$ be coprime integers with $0<e<d$, and let
$\zeta\in\CT$ be a primitive root of unity of degree $d$. Consider
an affine toric surface $ X_{d,\,e} \, = \, \AA^2/\, G_{d,\, e}$,
where $G_{d,\, e}=\langle g\rangle$ is the cyclic group generated
by an element $ g=\begin{pmatrix}
\zeta^e & 0 \\
0 & \zeta
\end{pmatrix}\in\GL(2,\CC)
$.  Let $Q=\pi(\bar 0)$ denote the unique singular point of
$X_{d,\,e}$, where $\pi \colon \AA^2 \to X_{d,\,e}$ is the
quotient morphism and  $\bar 0=(0,0)\in\AA^2$. We let
$N(G_{d,\,e})$ denote the normalizer of the subgroup $G_{d,\,e}$
in $\GL(2,\CC)$ and by $\cN(G_{d,\,e})$ that in the group $\Aut(
\AA^2)$.

\brem \label{rem1} By~\cite[Section~2.6]{Ful}, the surfaces
$X_{d,\,e}$ and $X_{d',\,e'}$ are isomorphic if and only if $d=d'$
and either $e=e'$ or $ee'=1 \mod d$. The latter two possibilities
are related via the twist $\tau:(x,y)\mapsto (y,x)$ on $\AA^2$.
\erem

For a linear form $l$  on $\AA^2$ and for $c\in\CC$ we let
$C_l(c)$ denote the affine line $l=c$.

\bthm \label{prop3}\bnum\item[(a)] Up to the action of the
automorphism group $\Aut (X_{d,\,e})$, any irreducible acyclic
curve $\Oc$ on $X_{d,\,e}$ which do not pass through the singular
point $Q\in X_{d,\,e}$ is equivalent either to $\pi(C_x(1))$ or to
$\pi(C_y(1))$. In particular, $\Oc$ coincides with an orbit of a
$\G_a$-action on $X_{d,\,e}$, and also with an orbit closure of a
$\G_m$-action  on $X_{d,\,e}$.
\item[(b)] The curves  $\pi(C_x(1))$ and $\pi(C_y(1))$ are
equivalent on $X_{d,\,e}$ if and only if $e^2 \equiv 1 \mod d$.
\enum \ethm

\bproof To show (a) we let $\Oc$ be an irreducible acyclic curve
on $X_{d,\,e}$ not passing through the singular point $Q$ and
$C=\pi^{-1}(\Oc)$ be its total preimage in $\AA^2$. The morphism
$\pi$ is finite of degree $d$, so any component $C^i$ of $C$ maps
to $\Oc$ properly. Since the cyclic group $G_{d,\, e}$ acts freely
on $\AA^2\setminus\{\bar 0\}$ and $C\subseteq\AA^2\setminus\{\bar
0\}$, the map $\pi|_{C^i}: C^i\to\Oc$ is a non-ramified covering.
However,  $\Oc$ being simply connected it does not admit any
non-trivial covering. Therefore $C$ has $d$ disjoint irreducible
components $C^1,\ldots,C^d$ mapped isomorphically onto $\Oc$ under
$\pi$. Furthermore, the cyclic group $G_{d,\, e}$ acts simply
transitively on the set $\{C^1,\ldots,C^d\}$.

Write a reduced defining equation of $C^1$ as $p-1=0$, where $p\in
\CC[x,y]$. Any regular invertible function on a connected and
simply connected variety is constant. Hence for every
$i=1,\ldots,d$ the restriction of $p|C^i$ is constant, say,
$\kappa_i$ i.e., $C^i\subseteq p^{-1}(\kappa_i)$. If $C^1$ were
singular then by theorem \ref{thm-LZ} it would be equivalent to a
curve $C_{a,b}=\{y^a-x^b=0\}$, where $\min\{a,b\}>1$. For $c\neq
0$ the Euler characteristic of the fiber $y^a-x^b=c$ is negative.
Hence this fiber cannot carry a curve with Euler characteristic
$1$. This leads to a contradiction, because $d>1$ by our
assumption. Thus the curve $C^1$ is smooth.

It follows that every fiber of $p$ is isomorphic to $\AA^1$. Hence
there is an automorphism $\delta\in\Aut(\AA^2)$  sending the
curves $C^i$ to the lines $y=\kappa_i$ with distinct $k_i$, where
$\kappa_1=1$. Moreover, we may suppose that $\delta(\bar 0)=\bar
0$. Letting $g'=\delta\circ g\circ \delta^{-1}$ we obtain $g'(\bar
0)=\bar 0$, ${g'}^d=g^d=\id$, and $g'(C_y)=C_y$ i.e.,
$g'\in\Stab(C_y)=\Jonq^+(\AA^2)$, see proposition \ref{prop1}.
Furthermore, $$\rho(g')=dg'(\bar 0)=d\delta(\bar 0)\circ g\circ
(d\delta(\bar 0))^{-1}\in\TT\,,$$ where
$\rho:\Jonq^+(\AA^2)\to\TT$ is the canonical surjection (see
proposition \ref{prop-quasitorus}). Hence the elements
$g=\rho(g)\in\TT$ and $\rho(g')\in\TT$ are conjugated in
$\GL(2,\CC)$ and so either $\rho(g')=\rho(g)=g$ or $\rho(g')=\tau
g \tau$.

Since $ \ord(g')=d>1$ we have $g'\notin U^+$. It follows from
lemma \ref{lem-uni-rad} and proposition \ref{prop-quasitorus} that
$\mu^{-1}g'\mu=\rho(g')$ for a suitable $\mu\in U^+$. In the case
where $\rho(g')=g$ we obtain $\mu^{-1}\delta g\delta^{-1}\mu=g$,
and in the case where $\rho(g')=\tau g\tau$ we get
$\tau\mu^{-1}\delta g\delta^{-1}\mu\tau=g$. In the former case
$\mu^{-1}\delta(C^1)=C_y(1)$, where
$\mu^{-1}\delta\in\cN(G_{d,e})$, while in the latter one
$\tau\mu^{-1}\delta(C^1)=C_x(1)$, where again
$\tau\mu^{-1}\delta\in\cN(G_{d,e})$. Since the corresponding
element normalizes the group $G_{d,\, e}$ it descends to an
automorphism of the surface $X_{d,\,e}$ which sends $\Oc$ to the
curve $\pi(C_y(1))$ in the former case and to $\pi(C_x(1))$ in the
latter one.

Furthermore, $\Oc$ is an orbit of a $\G_a$-action on $X_{d,\,e}$
induced by a $\G_a$-action on $\AA^2$ with $C^1$ as an orbit,
which commutes with the $G_{d,\,e}$-action and is defined via
$$t.(x,y)= (x, y+tx^{e'})\quad\text{if}\quad
\Oc=\pi(C_x(1))\quad\text{and}\quad t.(x,y) =(x+ty^e,
y)\quad\text{if}\quad \Oc=\pi(C_y(1))\,,$$ respectively, where
$t\in\G_a$. This shows  (a).

To show (b) we assume first that $e^2 \equiv 1 \mod d$ (and so
$e'=e$). Then the involution $\tau:(x,y)\mapsto (y,x)$ normalizes
the subgroup $G_{d,\, e}$ and induces an automorphism of
$X_{d,\,e}$ which interchanges the curves $\pi(C_x(1))$ and
$\pi(C_y(1))$ (and also $\pi(C_x)$ and $\pi(C_y)$).

Conversely, assume that there is an automorphism of $X_{d,\,e}$
which sends $\pi(C_x(1))$ to $\pi(C_y(1))$. It can be lifted to an
automorphism of the Cox ring $\CC[x,y]$ of the surface
$X_{d,\,e}$. Hence there is an element $\gamma \in \cN(G_{d,\,
e})$ which sends $C_x(1)$ to $C_y(1)$ and, moreover, sends the
variable $y$ to $x$. By proposition~\ref{prop1}(a) $\gamma$ has
the form $(x,y)\longmapsto (\alpha y + f(x),\,x)$. A direct
computation shows that
$$\gamma^{-1}\circ g \circ \gamma:
(x,y)\longmapsto (\zeta x + h(y),\,\zeta^e y)\,,$$ where
$\gamma^{-1}\circ g \circ \gamma\in G_{d,\, e}$ because $\gamma
\in \cN(G_{d,\, e})$. It follows that $h=0$ and
$(\zeta^e)^e=\zeta$ and so $e^2 \equiv 1 \mod d$. Now the proof is
completed. \eproof

\bcor\label{fam-curves} Every irreducible, acyclic curve $\Oc$ on
$X_{d,\,e}$ not passing through the singular point  $Q\in
X_{d,\,e}$ belongs to a pencil $\mathcal L$ consisting of
one-dimensional orbits of an effective $\G_a$-action on
$X_{d,\,e}$ and a fixed point curve $\cC_0$ passing through $Q$.
The members of  $\mathcal L$ different from $\cC_0$ are equivalent
under the $\G_m$-action on $X_{d,\,e}$ induced by the
$\TT_{1,1}$-action on $\AA^2$. \ecor

The union of several one-dimensional orbits of a $\G_a$-action on
$X_{d,\,e}$ is a disconnected, simply connected curve. In fact
every such curve arises in this way, as the reader can easily
derive from the previous results.

\bcor\label{disconn-curves} Every disconnected, simply connected
curve $\Oc$ on $X_{d,\,e}$ is equivalent under the
$\Aut(X_{d,\,e})$-action either to
$\bigcup_{i=1}^r\pi(C_x(\kappa_i))$, or to
$\bigcup_{i=1}^r\pi(C_y(\kappa_i))$, where $r\ge 2$ and
$\kappa_1,\ldots,\kappa_r\in\AA^1$ are distinct. \ecor


\subsection{Acyclic curves through the singular point}
\label{sec3}

In this subsection we describe the acyclic curves on a singular
affine toric surface $X_{d,e}$ ($d>1$) passing through the
singular point $Q$. We let as before $C_{a,\, b}=\{y^a-x^b=0\}$,
where $ a, b\ge 1$ and $\gcd(a,b)=1$. In particular
$C_{1,1}=\{x-y=0\}$. We keep the notation $C_x=\{x=0\}$  and
$C_y=\{y=0\}$.

Due to the following theorem, the set of all equivalence classes
of irreducible, acyclic curves on $X_{d,e}$ through $Q$ is
countable. A similar fact in the smooth case $X_{1,1}=\AA^2$ is
well known, see theorem \ref{thm-LZ}(a).

\bthm \label{prop4}\bnum\item[(a)] Up to the action of the
automorphism group $\Aut(X_{d,\,e})$, every irreducible acyclic
curve $\Oc$ on $X_{d,\,e}$ passing through the singular point
$Q\in X_{d,\,e}$ is equivalent to one of the curves $\pi(C_{a,\,
b})$, $\pi(C_x)$, or $\pi(C_y)$.\item[(b)] The curves  $\pi(C_x)$
and $\pi(C_y)$ are equivalent on $X_{d,\,e}$ if and only if $e^2
\equiv 1 \mod d$, if and only if the twist $\tau:(x,y)\longmapsto
(y,x)$ descends to an automorphism of the surface $X_{d,\,e}$.
\enum \ethm

\bproof (a) The curve $D=\pi^{*}(\Oc) \setminus \{\bar 0\}$ is
reduced and the projection $D \to \Oc \setminus \{Q\}$ is an
unramified cyclic covering of degree $d$. It follows that the
irreducible components $D^1, \ldots, D^r$ of $D$ are disjoint, and
for every $i=1,\ldots,r$ the restriction $\pi|_{D^i}:D^i\to \Oc
\setminus \{Q\}$ is an unramified cyclic covering. Furthermore,
the cyclic group $G_{d,\, e}$ acts transitively on the set of
these components i.e., the generator $g$ of $G_{d,\, e}$ permutes
them cyclically. In particular, $r$ divides $d$ and $D^i\simeq
\AT$ $\forall i$, while the closures $C^i=D^i\cup \{\bar 0\}$ and
$C=D\cup \{\bar 0\}=\pi^{-1}(\Oc)$ are acyclic. Clearly, the curve
$C$ cannot have more than one singular point, hence it cannot be
of type (III). According to the remaining types we distinguish the
following cases.

\smallskip

{\it Case 1:} \ {\em $r=1$ and $C=C^1$ is a smooth acyclic curve
of type (I).} Then $C^1$ is $G_{d,\, e}$-stable and passes through
the origin. Similarly as in the proof of theorem~\ref{prop3}, one
can show that in suitable new coordinates in $\AA^2$ we have
$C^1=C_y$ and $g$ acts diagonally either via $(x,y)\mapsto
(\zeta^e x, \zeta y)$, or via $(x,y)\mapsto (\zeta^{e'} x, \zeta
y)$. In the second case after transposition $(x,y)\mapsto (y,x)$
we obtain that $C^1=C_x$ and $g$ acts via $(x,y)\mapsto (\zeta x,
\zeta^{e'} y)$, that is by an element of the cyclic group
$G_{d,e}$. In any case up to the $\Aut(X_{d,\, e})$-action, the
curve $\Oc$ is equivalent either to $\pi(C_x)$ or to $\pi(C_y)$.

\smallskip

{\it Case 2\ :} \ {\em $r>1\,$ and $C$ is an acyclic curve of type
(II) or (IV) with an ordinary singularity at the origin.} By
theorem \ref{thm-LZ} a suitable automorphism
$\gamma\in\Aut(\AA^2)$ sends the  reduced curve $C=\pi^*(\Oc)$ to
a union $C'=C'^1+\ldots+ C'^r$ of affine lines through the origin
given by equation \be\label{can-form-lin}
y(y-\kappa_2x)\ldots(y-\kappa_{r}x)=0, \quad \text{where} \quad
\kappa_i\in\CC^\times\quad \text{are distinct}\,. \ee The curve
$C'$ is stable under the action  on $\AA^2$ of the cyclic group
$\gamma G_{d,e} \gamma^{-1}=\langle g'\rangle$, where $g'=\gamma g
\gamma^{-1}\in\Aut(\AA^2)$. By corollary \ref{cross-lin}
$\Stab(C')\subseteq\GL(2,\CC)$, hence $g'\in\GL(2,\CC)$.

There exists an element $\delta\in\GL(2,\CC)$ such that
 $\delta g' \delta^{-1}=g$ is diagonal and acts via
$(x,y)\mapsto (\zeta^e x, \zeta y)$. Since no component $C^i$ of
$C$ is stable under $g$, the composition $\delta\gamma$ sends each
$C^i$ to a line through the origin different from a coordinate
axis. Since all such lines are $\TT$-equivalent, their images in
the surface $X_{d,e}$ are also equivalent under the action on
$X_{d,e}$ of the quotient torus $\TT'=\TT/G_{d,e}$.

The resulting automorphism $\delta\gamma$ from the centralizer of
the subgroup $G_{d,e}$ in $\Aut(\AA^2)$ rectifies $C$ and sends
$C^1$ to a line $\TT$-equivalent to $C_{1,1}$. Consequently, the
curve $\Oc$ on $X_{d,e}$ is equivalent to $\pi(C_{1,\, 1})$ under
the $\TT'$-action on $X_{d,e}$ and the automorphism
$\pi_*(\delta\gamma)\in\Aut(X_{d,e})$.

\smallskip

{\it Case 3:}  \ {\em $r>1$ and $C$ is an acyclic curve of type
(V) with smooth components $C^i$ and a non-ordinary singularity at
the origin.} By theorem \ref{thm-LZ} in this case a suitable
automorphism $\gamma\in\Aut(\AA^2)$ sends the reduced curve
$C=C^1+\ldots+ C^r$ on $\AA^2$ to a curve $C'$ given by equation
\be\label{eq5} x^{\varepsilon_x}y(y-\kappa_2x^b)
\ldots(y-\kappa_{r}x^b)=0\quad\text{with distinct}\quad
\kappa_i\in\CC^\times\,,\ee where $b>1$ and
$\varepsilon_x\in\{0,1\}$. The group $G'_{d,e}=\gamma
G_{d,e}\gamma^{-1}\subseteq \Aut(\AA^2)$ acts transitively on the
set of components $C'^i$ of $C'$. Since the $G'_{d,e}$-action
preserves tangency we have $\varepsilon_x=0$. We may also suppose
that $k_2=1$. Since $G'_{d,e}\subseteq\Stab(C')$, the singular
point $\bar 0\in C'$ is fixed under the $G'_{d,e}$-action on
$\AA^2$.

There is an element $g'\in G'_{d,e}$ which sends the component
$C'^2=\{y-x^b=0\}$ to $C'^1=C_y$. Since $h:(x,y)\mapsto (x,y-x^b)$
does the same, according to proposition \ref{prop1}(a) $g'$ can be
written as
$$g':(x,y)\longmapsto (\alpha x+f(y-x^b),\, \beta(y-x^b))\,$$
for some $\alpha,\beta\in\CC^\times$ and $f\in\CC[z]$. Hence $g'$
maps the affine line $C'^1=C_y$ to the parameterized curve
$\{(\alpha t+f(-t^b),-\beta t^b)\,|\,t\in\AA^1\}$. However
$g'(C'^1)=C'^i$ for some $i\in\{2,\ldots,r\}$. It follows that
$f=0$. Therefore $g':(x,y)\mapsto (\alpha x, \beta(y-x^b))$ and so
$$g'^k:(x,y)\mapsto \left(\alpha^k x, \beta^k y
-\beta x^b(\beta^{k-1}+\alpha^b\beta^{k-2}
+\ldots+\alpha^{b(k-1)}) \right)\,.$$ Since $g'^k=\id$ for some
$k|d$, we have $\alpha^k=\beta^k=1$ and
$(\beta/\alpha^b)^{k-1}+\ldots+1=0$. This implies that
$\beta\neq\alpha^b$.

Thus the triangular automorphism
$$h':(x,y)\mapsto \left(x,y+\frac{\beta}{\alpha^b-\beta} x^b\right)\,$$
is well defined and sends the curve $C'$ to a new one $C''=\sum_i
C''^i$ given by a similar equation. Furthermore, $h'$ conjugates
$g'$ with
$$g'': (x,y)\longmapsto (\alpha x,
\beta y)\,,$$ where $(\alpha,\beta)$ can be written either as
$(\xi^e,\xi)$ or as $(\xi^{e'},\xi)$ for a primitive $d$th root of
unity $\xi$ (recall that $ee'\equiv 1\mod d$). In the latter case
we apply additionally the transposition of coordinates to get
$(\alpha,\beta)=(\xi,\xi^{e'})$. The composition $h'\gamma$
normalizes the group $G_{d,e}$ in $\Aut(\AA^2)$ and sends $C^1$ to
a member of the pencil $\mathcal L$, where as before  $$\mathcal
L=\{y-\kappa x^b=0\}_{\kappa\in\CC}\cup\{bC_x\}\,.$$ However,
every member of $\mathcal L$ different from the coordinate axes
$C_x,\,C_y$ is equivalent to the curve $C_{1,b}$ under the
$\TT$-action on $\AA^2$. Finally $h'\gamma$ induces an
automorphism $\theta$ of the quotient surface $X_{d,e}$ such that
the image under $\theta$ of the curve $\Oc=\pi(C^1)$ is
$\TT'$-equivalent either to $\pi(C_y)$ or to $\pi(C_{1,b})$, as
required.

\smallskip

{\it Case 4:} \ {\em $C$ is an acyclic curve of type (VI) with all
the components $C^i$ being singular.} By theorem \ref{thm-LZ} a
suitable automorphism $\gamma\in\Aut(\AA^2)$ sends $C$ to a curve
$C'$ given by equation (\ref{eq-par2}), where
$\varepsilon_x=\varepsilon_y=0$, $a,b>1$, and $\gcd(a,b)=1$. We
may also assume that $\kappa_1=1$ and so $C'^1=C_{a,\, b}$. Let as
before $g'=\gamma g \gamma^{-1}$. Since the $\TT$-action is
transitive on the members of the pencil
$$\mathcal L=\{y^a-\kappa x^b=0\}_{\kappa\in\CC}\cup\{bC_x\}$$
different from the coordinate axes $C_x$ and $C_y$, there is an
element $h\in\TT$ such that $g''=h\circ g'$ preserves the curve
$C_{a,\, b}$. By proposition~\ref{prop2} $g''\in\TT$, hence also
$g'\in\TT$. Since the eigenvalues $\zeta^e$ and $\zeta$ of $g$ and
$g'$ are the same, up to interchanging the coordinates we have
$g'\in G_{d,e}$. Hence $\gamma$ normalizes the group $G_{d,e}$.
Reasoning as before we conclude that the curve $\Oc$ on the
surface $X_{d,e}$ is equivalent to $\pi(C_{a,\, b})$, as stated.
This shows (a). The proof of (b) goes in the same way as that of
theorem \ref{prop3}(b) and so we leave it to the reader. \eproof

\brems\label{cano} 1. If $e>1$ then the curve
$\pi^{-1}(\pi(C_x))=C_x$ on $\AA^2$ is irreducible, while
$\pi^{-1}(\pi(C_{1,\,1}))$ is reducible. Hence these plane curves
are not equivalent and so their images $\pi(C_x)$ and $\pi(C_{1,\,
1})$ on the surface $X_{d,\, e}$ are not equivalent either (under
the $\Aut(X_{d,\, e})$-action).

2. The following simple example shows that, in contrast with
theorem \ref{prop3}, a curve on $X_{d,e}$ isomorphic to the affine
line and passing through the singular point can be non-equivalent
to the image of one of the curves $C_{1,1},\,C_x$, and $C_y$.

Indeed, the curve $C_{2,\, 3}$ in $\AA^2$ is singular.
Nonetheless, its image $\cC$ in the affine toric surface $X_{5,\,
4}$ passes through the singular point and is isomorphic to the
affine line $\AA^1$. Clearly, $\cC$ is not equivalent in $X_{5,\,
4}$ to $\pi(C_{1,1})$, $\pi(C_x)$, or $\pi(C_y)$. \erems


\subsection{Acyclic curves as
orbit closures} \label{sec4}

Summarizing the  results of the previous subsections we arrive at
the following alternative description.

\bthm \label{tmain} Let $X$ be an affine toric surface over $\CC$
with the acting torus $T$. Then every irreducible acyclic curve
$\Oc$ on $X$ coincides with the closure of a non-closed orbit of a
regular $\G_m$-action on $X$. \footnote{Clearly, such an orbit
closure is acyclic.} Furthermore, up to an automorphism of $X$,
such a curve $\Oc$ is the closure of a non-closed orbit of a
subtorus of $T$. \ethm

\bproof Let $\Oc$ be an irreducible  acyclic curve on a toric
surface $X$. If $X$ is smooth then this is one of the surfaces
$\AT\times \AT$, $\AA^1\times \AT$,  or $\AA^2$. There is no
acyclic curve on $X=\AT\times \AT$, and the only irreducible
acyclic curves on $\AA^1\times \AT$ are of the form $\AA^1 \times
\{\text{pt}\}$. Hence our assertion holds for these surfaces. In
the case of the affine plane $X\simeq\AA^2$ the result follows
from theorems \ref{thm-AMS} and \ref{thm-LZ}.

By~\cite[\S2.2]{Ful}, every singular affine toric surface is
isomorphic to one of the surfaces $X_{d,\, e}=\AA^2/\, G_{d,\,e}$,
where $d>1$. Finally in this case the result follows from
theorems~\ref{prop3} and \ref{prop4}. \eproof

\bcor\label{co1} Any irreducible acyclic curve on an affine toric
surface has at most one singular point. If the surface is singular
then this point coincides with the singular point of the surface.
\ecor

The following result generalizes theorem \ref{thm-AMS} of
Abhyankar, Moh, and Suzuki.

\bthm\label{co2} Up to the action of the group $\Aut(X_{d,e})$
there are only finitely many different embeddings
$\AA^1\hookrightarrow X_{d,e}$. \ethm

\bproof It suffices to show that the smooth curves  on $X_{d,\,
e}$ of the form $\pi(C_{a,\, b})$ belong to a finite set of
equivalence classes. Notice that the subalgebra $\CC[x,y]^{G_{d,\,
e}}$ of $G_{d,\, e}$-invariants is generated by the monomials
$y^d, xy^{c_1},\ldots,x^{d-1}y^{c_{d-1}},x^d$, where $0<c_k<d$ and
$c_k+ke\equiv 0\mod d$. These monomials define a closed embedding
$X_{d,\, e}\hookrightarrow \AA^{d+1}$. The image of the curve
$\pi(C_{a,\, b})$ under this embedding is
$$
(t^{db}, t^{a+c_1b},\ldots, t^{(d-1)a+c_{d-1}b}, t^{da}), \quad
t\in\CC\,.
$$
This image is smooth if and only if one of the exponents, say,
$\delta$ of our monomials coincides with the greatest common
divisor of all the exponents. Since $a$ and $b$ are coprime
$\delta|d$.

In the case where $\delta=ka+c_kb$ for some $k\ge 1$ we obtain
$\delta=ka+c_kb\le d$ and so $a+b\le d$. The number of all
possible such pairs $(a,b)$ is finite.

If $\delta\in\{da,db\}$ then $a=1$ or $b=1$ because $\gcd(a,b)=1$.
Suppose for instance that $a=1$ and $\delta=d$, the other case
being similar. Thus
$$ek+c_k\equiv 0\mod d\quad\text{and}\quad k+bc_k\equiv 0\mod d
\,$$ for all possible values of $k\ge 1$. For $k=1$ it follows
that $eb\equiv 1\mod d$. Hence the curve $C_{a,b}=C_{1,b}$ is
stable under the $G_{d,e}$-action on $\AA^2$. The automorphism
$$(x,y)\longmapsto (x,y-x^b)$$ commutes with the $G_{d,e}$-action
and sends this curve to the axis $C_y$. Therefore for any $b\ge 1$
the curves $\pi(C_{1,b})$ and $\pi(C_y)$ on $X_{d,e}$ are
equivalent. Now the proof is completed. \eproof

\brem\label{DBH} Consider an affine toric variety $X$ of dimension
$n$. It is known (see~\cite{BH}) that if an $(n-1)$-dimensional
torus~$T'$ acts effectively on $X$, then $T'$ is conjugate in the
group $\Aut(X)$ to a subtorus of the acting torus $T$ of $X$. An
analogous result for tori of codimension $\ge 2$ is unknown. We
conclude section \ref{sec-acyc-curves} by the following related
problem. \erem

\medskip

\noindent {\bf Questions.} Let $X$  be an affine toric variety of
dimension $n$ with acting torus $T$. Assume that a closed
subvariety $Y$ of $X$ admits in turn an action of an algebraic
torus $T'$ with an open orbit.

\smallskip

\noindent {\em Is it true that $Y$ can be realized as the orbit
closure of a $T'$-action on $X$? Is, moreover, this $T'$-action on
$X$ conjugated to the action on $X$ of a subtorus of the acting
torus $T$?}

\smallskip

\noindent This is indeed the case for $n=2$ as follows from
theorem~\ref{tmain}.


\subsection{Reducible acyclic curves
on affine toric surfaces}\label{reducible} Let us start with the
following lemma.

\blem\label{simply-conn-cover} Consider the affine toric surface
$X=X_{d,e}$ with the quotient map $\pi:\AA^2\to X=\AA^2/G_{d,e}$.
If $\Oc$ is a reduced, simply connected curve on $X$, then the
total transform $C=\pi^*(\Oc)$ of $\Oc$ in $\AA^2$ is also reduced
and simply connected. \elem

\bproof If $d=1$ i.e., $X\simeq\AA^2$, or the curve $\Oc$ is
irreducible, then the assertion follows by the same argument as in
the proof of theorems \ref{prop3} and \ref{prop4}. Assume further
that $d>1$ and $\Oc$ is reducible. Letting $Q=\pi({\bar 0})\in X$
and $\pi^*(\Oc)=C=C^1+\ldots+C^s$, where every irreducible
component $C^i$ of $C$ is simply connected, we consider the
following cases.

\smallskip

\noindent {\em Case 1 : Every irreducible component $\Oc^k$ of
$\Oc$ passes through $Q$.} Then any two such components meet only
at $Q$, and it is easily seen that any two distinct components
$C^i$ and $C^j$ of $C$ also meet only at the origin. In this case
$C$ is connected and simply connected i.e., acyclic.

\smallskip

\noindent {\em Case 2 : There are two distinct crossing
components, say, $\Oc^k$ and $\Oc^l$ of $\Oc$ not passing through
$Q$.} According to theorem~\ref{prop3} (cf. also Corollaries
\ref{fam-curves} and \ref{disconn-curves}), under the action of
the normalizer $\On(G_{d,e})$ on $\AA^2$ the total transform
$\pi^*(\Oc^k)$ (or $\pi^*(\Oc^l)$) is equivalent to a union of $d$
parallel lines, which are parallel either to $C_x$ or to $C_y$ and
are cyclically permuted under the $G_{d,e}$-action. Furthermore,
each component of $\pi^*(\Oc^k)$ meets every component of
$\pi^*(\Oc^l)$ in $d$ distinct points. All these $d^2$
intersection points must project to the unique intersection point
$\Oc^k\cap\Oc^l$. Hence they should belong to a $G_{d,e}$-orbit,
which is impossible. This contradiction shows that the components
of $\Oc$ not passing through $Q$ do not meet.

\smallskip

\noindent {\em Case 3 : There is just one component, say, $\Oc^1$
of $\Oc$ not passing through $Q$.} The union of the other
components of $\Oc$ can meet $\Oc^1$ in at most one point, and
they meet each other at $Q$. Hence either $\Oc^1$ does not meet
this union, or there is just one component, say, $\Oc^2$ of $\Oc$
passing through $Q$ which meets $\Oc^1$. In the former case the
reduced curve $C=\pi^*(\Oc)$ is clearly simply connected, as
stated. In the latter case we may assume as before that
$\pi^*(\Oc^1)$ is a union of $d$ lines parallel to a coordinate
axis and cyclically permuted under the $\G_{d,e}$-action. Every
component of the total transform $\pi^*(\Oc^2)$ is simply
connected, passes through the origin, and meets one of these
lines. Hence it meets all the parallel lines. If the curve
$\pi^*(\Oc^2)$ is reducible then its irreducible components meet
only at the origin, and meet one of the lines in  $\pi^*(\Oc^1)$
in at least two distinct points. These points project in $X$ to
distinct smooth points. The latter contradicts the assumption of
simply connectedness of the curve $\Oc$, because in this case we
obtain a non-contractible cycle in $\Oc$. Hence the curve
$\pi^*(\Oc^2)$ is irreducible and meets every line in
$\pi^*(\Oc^1)$  in just one point, while the other curves
$\pi^*(\Oc^k)$, $k\ge 3$, do not meet these lines at all. This
shows that $C$ is simply connected.

\noindent {\em Case 4 : There are two disjoint components, say,
$\Oc^1$ and $\Oc^2$ of $\Oc$ not passing through $Q$.} Since the
total transforms $\pi^*(\Oc^1)$ and $\pi^*(\Oc^2)$ are disjoint,
they can be simultaneously transformed into unions of lines
parallel to the same coordinate axis. Hence every component, say,
$\Oc^3$ of $\Oc$ passing through $Q$ and meeting $\Oc^1$ meets
also $\Oc^2$, and vice versa.  Using the same argument as before
it is easily seen that there could be at most one such component
$\Oc^3$, and the total preimage $C=\pi^*(\Oc)$ is simply
connected, as required.

This ends the proof. \eproof

Using this lemma, in the following theorem we give a description
of all reduced, simply connected curves on affine toric surfaces.

\bthm\label{thm-reducible-curves} Every reduced, acyclic curve
$\Oc$ on $X=X_{d,e}$ is equivalent to a curve $\pi(C)$, where
$C\subseteq \AA^2$ is given by one of the equations \be\label{tri}
x\prod_{i=1}^r(y^a-\kappa_ix^b)=0\quad\text{or} \quad
y\prod_{i=1}^r(x^a-\kappa_iy^b)=0,\quad\text{where} \quad b\ge
0,\,\,\, a\ge 1,
 \quad\gcd(a,b)=1\,,\ee and  where $\kappa_i\in\CC$
 ($i=1,\ldots,r$) are pairwise distinct.
\ethm

\bproof Indeed, by lemma \ref{simply-conn-cover} the reduced plane
curve $C=\pi^*(\Oc)$ is simply connected. However, by theorem
\ref{thm-LZ}  every reduced, acyclic curve in $\AA^2$ is given in
appropriate coordinates by equation (\ref{tri}), while every
reduced, disconnected, simply connected affine plane curve is
equivalent to a finite union of parallel lines. It remains to show
that, up to a permutation of the symbols $x$ and $y$, the
corresponding coordinate change can be chosen in the normalizer
$\On(G_{d,e})$ of the group $G_{d,e}$ in $\Aut(\AA^2)$. The latter
can be done in the same way as in the proof of
theorems~\ref{prop3} and \ref{prop4}. We leave the details to the
reader. \eproof


\section{Automorphism groups of affine toric surfaces}
\label{sec-gr-aut} In this section we prove an analog of the
Jung-van der Kulk theorem \ref{thm-JvdK} for affine toric surfaces
and study algebraic groups acting on such surfaces.

\subsection{Free amalgamated product structure}\la{next}
Consider again an affine toric surface $X_{d,e}=\AA^2/G_{d,e}$. We
assume as usual that $1\le e<d$, $\gcd(d,e)=1$, and
$G_{d,e}=\langle g\rangle$, where
$$g:(x,y)\longmapsto (\zeta^ex,\zeta y)\quad\text{with}
\quad \zeta=\exp\left(\frac{2\pi i}{d}\right)\,.$$

\bnota\la{nota-norm} Let $G\subseteq \GL(2,\CC)$. Letting as
before $N(G)$ denote the normalizer of $G$ in the group
$\GL(2,\CC)$ and $\cN(G)$ that in the group $\Aut(\AA^2)$, we
abbreviate $$N_{d,e}=N(G_{d,e})\text{ and }
\cN_{d,e}=\cN(G_{d,e})\,.$$ It is easily seen that
\[
N_{d,e}=\begin{cases} \GL(2,\CC) & \text{ if   $e=1$,}\\
N(\TT)=\langle\TT,\tau\rangle  &\text{ if  $ e>1$  and
$e^2\equiv 1\mod d$,}\\
\TT&\text{ otherwise, } \end{cases}
\]
where $\tau:(x,y)\longmapsto (y,x)$ is a twist and $\TT$ stands
for the maximal torus in $\GL(2,\CC)$ consisting of the diagonal
matrices. We let $B^\pm$ denote the Borel subgroup of all upper
(lower, respectively) triangular matrices in $\GL(2,\CC)$.
Consider the  subgroups \be\label{nde} N_{d,e}^\pm=N_{d,e}\cap
{\Jonq}^{\pm}(\AA^2)=\begin{cases}
B^\pm &\text{ if $e=1$,}\\
\TT&\text{ otherwise}
\end{cases}
\ee and \[ \cN_{d,e}^\pm=\cN_{d,e}\cap {\Jonq}^{\pm}(\AA^2)\,.\]
The latter subgroups are described in lemma
\ref{lem-descr-normalis-bis} below. Notice that \be\label{tor-int}
\cN_{d,e}^+\cap\cN_{d,e}^-= {\Jonq}^{+}(\AA^2) \cap
{\Jonq}^{-}(\AA^2)=\TT\,.\ee \enota

With this notation we can state an analog of theorem
\ref{thm-JvdK} by Jung and van der Kulk.

\bthm\la{thm-JvdK-toric-surf} If $e^2\not\equiv 1 \mod d$ then
\be\la{eq-amalg-1} \Aut (X_{d,e}) \simeq
\cN_{d,e}^+/G_{d,e}*_{\TT/G_{d,e}} \cN_{d,e}^-/G_{d,e}\,, \ee
while for $e^2\equiv 1 \mod d$ we have \be\la{eq-amalg-2} \Aut
(X_{d,e}) \simeq \cN_{d,e}^+/G_{d,e}*_{N_{d,e}^+/G_{d,e}}
N_{d,e}/G_{d,e}\,. \ee \ethm

There should be possible to derive this theorem by using the
techniques elaborated by Danilov and Gizatullin \cite{DG}.
However, we prefer a direct approach through an equivariant
version \ref{lem-AMS2} of the Abhyankar-Moh-Suzuki theorem. Within
this approach theorem \ref{thm-JvdK-toric-surf} is an immediate
consequence of lemma \ref{lem-nornalis} and proposition
\ref{prop-amalga} below.

\blem\la{lem-nornalis} There is an isomorphism
$$ \Aut (X_{d,e}) \simeq \cN_{d,e}/G_{d,e}\,.$$
\elem

\bproof The affine plane $\AA^2$ can be viewed as the spectrum of
a Cox ring of the toric surface $X_{d,e}$, see \cite{Cox} or
\cite[I.6.1]{ADHL}. Hence every automorphism
$\Phi\in\Aut(X_{d,e})$ can be lifted (in a non-unique way) to an
element $\varphi\in\cN_{d,e}$. \footnote{Alternatively, this
follows from the monodromy theorem, see the proof of theorem
\ref{thm-non-abel-aut} below.} This yields an exact sequence (see
\cite[Thm. 5.1]{AG})
\[1\to G_{d,e}\to\cN_{d,e}\to \Aut(X_{d,e})\to 1\,,\] as claimed.
\eproof

\bprop\la{prop-amalga} If $e^2\not\equiv 1 \mod d$ then
\be\la{eq-amalg-1bis}  \cN_{d,e} \simeq \cN_{d,e}^+*_{\TT}
\cN_{d,e}^-\,, \ee while for $e^2\equiv 1 \mod d$
\be\la{eq-amalg-2bis} \cN_{d,e}  \simeq {\cN_{d,e}^+}*_{N_{d,e}^+}
N_{d,e} \,. \ee \eprop

The proof is done in lemmas
\ref{lem-descr-normalis-bis}--\ref{lem-ge-2} below.

\smallskip

Theorem \ref{thm-JvdK-toric-surf} follows now from lemma
\ref{lem-nornalis} and proposition \ref{prop-amalga} due to the
fact that the subgroup $G_{d,e}$ is normal in every group that
participates in (\ref{eq-amalg-1bis}) and (\ref{eq-amalg-2bis}).
Indeed, this can be seen directly or, alternatively, derived as a
simple consequence of a theorem by Cohen \cite{Coh} on
preservation of the  free amalgamated product structure in the
quotient.

\medskip

Recall (see remark \ref{rem-grading}) that the polynomial ring
$A=\CC[t]$ possesses a $\ZZ/d\ZZ$-grading
$$A=\bigoplus_{i=0}^{d-1} A_{d,i},
\quad\text{where}\quad A_{d,i}=t^i\CC[t^d]\,.$$ In terms of this
grading the normalizer $\cN_{d,e}^\pm$ admits the following
description.

\blem\la{lem-descr-normalis-bis} \bnum\item[(a)] The group
$\cN_{d,e}^+$ ($\cN_{d,e}^-$, respectively) consists of all de
Jonqi\`eres transformations $\varphi^+$ as in (\ref{eq-Jonq-plus})
($\varphi^-$ as in (\ref{eq-Jonq-minus}), respectively) with $f\in
A_{d,e}$ ($f\in A_{d,e'}$, respectively).
\item[(b)] The subgroup $\cN_{d,e}^\pm$
is the centralizer of $G_{d,e}$ in the group $\Jonq^\pm(\AA^2)$.
\enum \elem

\bproof We stick to the plus-case, the proof in the other one
being similar. We have $$\varphi^+\circ g\circ (\varphi^+)^{-1}:
(x,y)\longmapsto \left(\zeta^ex+f(\zeta\frac{y}{\beta}) -\zeta^e
f(\frac{y}{\beta}), \zeta y\right)\,.$$ Hence $\varphi^+\circ
g\circ  (\varphi^+)^{-1}\in G_{d,e}$ if and only if $f(\zeta
t)=\zeta^ef(t)$, if and only if $f\in A_{d,e}$. This shows (a). In
the latter case $\varphi^+\circ g\circ  (\varphi^+)^{-1}=g$, so
(b) follows. \eproof

For a pair of polynomials $\varphi=(u,v)$ we let $\deg
\varphi=\max\{\deg u,\,\deg v\}\,.$ The following result is an
immediate consequence of lemma \ref{lem-descr-normalis-bis}.

\blem\la{lem-descr-normalis-bis-bis} Assume as before that $1\le
e<d$ and $\gcd (d,e)=1$. Then the following hold. \bnum\item[(a)]
$\varphi(\bar 0)=\bar 0$ $\forall \varphi\in \cN_{d,e}^\pm$.
\item[(b)] $\cN_{d,e}^\pm\cap\Aff(\AA^2)=N_{d,e}^\pm=
\begin{cases} B^\pm, & e=1,\\
\TT, & e>1.\end{cases}$
\item[(c)] Let $\varphi^\pm$ be as in (\ref{eq-Jonq-plus}) and
(\ref{eq-Jonq-minus}), respectively. Assume that
$\varphi^\pm\not\in\TT$. Then $\deg \varphi^+\ge e$ and $\deg
\varphi^-\ge e'$, where $1\le e'<d$ and $ee'\equiv 1\mod d$. \enum
\elem

The next two lemmas provide a $G_{d,e}$-equivariant version of the
Abhyankar-Moh-Suzuki theorem \ref{thm-AMS}.

\blem\la{lem-AMS1} Let $C$ be a smooth, polynomial curve in
$\AA^2$ parameterized via $t\longmapsto (u(t),\,v(t))$, where
$u,v\in t\CC[t]$. If $G_{d,e}\subseteq\Stab(C)$ then either
$(u,v)\in A_{d,e}\times A_{d,1}$ or $(u,v)\in A_{d,1}\times
A_{d,e'}$. \elem

\bproof The tangent vector  to $C$ at the origin is
$w=(u'(0),\,v'(0))\in\AA^2$. Since it is stable under the tangent
$G_{d,e}$-action then either $g(w)=\zeta^e w$ or $g(w)=\zeta w$.
Thus  $g|C:C\to C$ acts either via $t\longmapsto \zeta^e t$ or via
$t\longmapsto \zeta t$. In the former case
$$g\circ(u,v)(t)=(\zeta^e u(t),\,\zeta v(t))=(u(\zeta^e t),v(\zeta^e t))\,,$$
and in the latter one
$$g\circ(u,v)(t)=(\zeta^e u(t),\,\zeta v(t))=(u(\zeta t),v(\zeta t))\,.$$
Now the assertion follows. \eproof

\blem\la{lem-AMS2} For a curve $C$ as in lemma \ref{lem-AMS1}
there is an automorphism
$\varphi\in\langle\cN_{d,e}^+,\,\cN_{d,e}^-\rangle$ which sends
$C$ to one of the coordinate axes $C_x$ and $C_y$. \elem

\bproof If $u=0$ or $v=0$ there is nothing to prove. Thus we may
suppose that $\deg u\ge \deg v>0$. By the Abhyankar-Moh-Suzuki
theorem \ref{thm-AMS} we have $\deg u=n\deg v$ for some $n\in\NN$.
By virtue of lemma \ref{lem-AMS1} either $(\deg u, \deg v)\equiv
(e,1)\mod d$ or $(\deg u, \deg v)\equiv (1,e')\mod d$. In both
cases it follows that $n\equiv e\mod d$ and so
$\varphi_1\in\cN_{d,e}^+$, where  $\varphi_1:(x,y)\longmapsto
(x-cy^n,y)$ (see  lemma \ref{lem-descr-normalis-bis}(a)). We have
$\varphi_1(u,v)=(u_1,v_1)=(u-cv^n,v)$. So we can choose
$c\in\CC^\times$ in such a way that $\deg u_1<\deg u$. We can
continue this procedure recursively until we reach one of the
pairs $(u_s,v_s)=(\alpha t, 0)$ or $(u_s,v_s)=(0,\beta t)$, where
$\alpha,\beta\in\CC^\times$. Then the product
$\varphi=\varphi_s\circ\ldots\circ\varphi_1$ is a required
automorphism. \eproof

\blem\la{lem-stabil-axes} For any $\varphi\in\cN_{d,e}$ we have
$\varphi(\bar 0)=\bar 0$ and \be\la{eq-inclu}
G_{d,e}\subseteq\Stab(\varphi(C_x)) \cap\Stab(\varphi(C_y))\,.\ee
In particular, \be\la{eq-equa} \cN_{d,e}\cap\Aff(\AA^2)
=N_{d,e}\,.\ee \elem

\bproof Since $\varphi$ normalizes the subgroup $G_{d,e}=\langle
g\rangle$ we have $\varphi^{-1}\circ g\circ\varphi=g^k$ for some
$k\in\NN$ and so $g\circ\varphi=\varphi\circ g^k$. Hence
$g(\varphi(C_x))=\varphi(g^k(C_x))=\varphi(C_x)$ and, similarly,
$g(\varphi(C_y))=\varphi(C_y)$. This yields (\ref{eq-inclu}). From
$$g(\varphi(\bar 0))=\varphi(g^k(\bar 0))=\varphi(\bar 0)$$ we
deduce that $\varphi(\bar 0)\in(\AA^2)^g=\{\bar 0\}$ i.e.
$\varphi(\bar 0)=\bar 0$. Now the last assertion follows easily.
\eproof

\blem\la{lem-generat} If $e^2\not\equiv 1\mod d$ then
\be\la{eq-plus-minus} \cN_{d,e}=\langle
\cN_{d,e}^+,\,\cN_{d,e}^-\rangle\,,\ee while for $e^2\equiv 1\mod
d$, \be\la{eq-plus-tau} \cN_{d,e}=\langle
\cN_{d,e}^+,\,N_{d,e}\rangle=\langle
\cN_{d,e}^+,\,\tau\rangle\,.\ee \elem

\bproof For $\varphi\in\cN_{d,e}$ we let $C=\varphi^{-1}(C_y)$. By
lemma \ref{lem-stabil-axes} (see (\ref{eq-inclu})) the cyclic
group $G_{d,e}$ stabilizes $C$ and $\bar 0\in C$. By virtue of
lemma \ref{lem-AMS2} there is an automorphism
$\psi\in\langle\cN_{d,e}^+,\,\cN_{d,e}^-\rangle$ which sends $C$
to one of the coordinate axes $C_x,\,C_y$. Letting
$\gamma=\psi\circ\varphi^{-1}\in\cN_{d,e}$ we get
$\gamma(C_y)=\psi(C)\in\{C_x,\,C_y\}$.  If $\gamma(C_y)=C_x$ then
the images of the coordinate axes $\pi(C_x)$ and $\pi(C_y)$ are
equivalent in the surface $X_{d,e}$. According to theorem
\ref{prop4} in this case $e^2\equiv 1\mod d$ and $\tau\in
N_{d,e}$.

Thus $\gamma(C_y)=C_y$ (i.e. $\gamma\in\Stab(C_y)$) if
$e^2\not\equiv 1\mod d$. By proposition \ref{prop1} in this case
$$\Stab(C_y)\cap\cN_{d,e}={\Jonq}^+(\AA^2)\cap\cN_{d,e}
=\cN_{d,e}^+\,.$$ Hence $\varphi\in \langle
\cN_{d,e}^+,\,\cN_{d,e}^-\rangle$ and so (\ref{eq-plus-minus})
follows.

Assume further that $e^2\equiv 1\mod d$ and $\gamma(C_y)=C_x$.
Then $\tau\circ\gamma(C_y)=C_y$ and so
$\tau\circ\gamma\in\cN_{d,e}^+$ and $\gamma\in\langle
\cN_{d,e}^+,\tau\rangle$. It follows that
$$\varphi=\gamma^{-1}\circ\psi\in\langle \cN_{d,e}^+,\,
\cN_{d,e}^-,\,\tau\rangle=\langle
\cN_{d,e}^+,\,\tau\rangle=\langle
\cN_{d,e}^+,\,N_{d,e}\rangle\,.$$ Now the proof is completed.
\eproof

We need the following analog of lemma 4.1 in \cite{Ka2}, see also
\cite[Theorem 5.3.1]{Wr1} and \cite[Lemma 1.9]{Wr3}. For the
reader's convenience we provide a short argument.

\blem\la{lem-generat-bis} Let as before $1\le e<d$, where $\gcd
(d,e)=1$. Consider an automorphism $\varphi\in\Aut(\AA^2)$ with
components $u,v\in\CC[x,y]$, written as an alternating product
\be\la{eq-alt-pro}
\varphi=\varphi_s\cdot\ldots\cdot\varphi_1\quad\text{with}\quad
\varphi_i\in\cN_{d,e}^\pm\setminus N^\pm_{d,e}\quad\text{and}\quad
\varphi_{i+1}\in\cN_{d,e}^\mp\setminus N^\mp_{d,e},\quad
i=1,\ldots,s-1\,,\ee where $s\ge 1$. Then
$$\deg u> \deg v\quad\text{ if}\quad
\varphi_s\in\cN_{d,e}^+\quad\text{and}\quad \deg u< \deg
v\quad\text{ if}\quad \varphi_s\in\cN_{d,e}^-\,.$$ In both cases
\be\la{eq-deg-prod} \deg\varphi:=\max\{\deg u,\deg
v\}=\prod_{i=1}^s\deg\varphi_i\,.\ee \elem

\bproof Both assertions are evidently true if $s=1$. Letting $s>1$
we assume by induction that  they hold for the product
$\psi=\varphi_{s-1}\cdot\ldots\cdot\varphi_1=(\tilde u, \tilde
v)$. Thus \be\la{eq-induct} \deg\psi
=\prod_{i=1}^{s-1}\deg\varphi_i\,,\ee $\deg \tilde u< \deg \tilde
v$  if $\varphi_{s-1}\in\cN_{d,e}^-$ (i.e.
$\varphi_{s}\in\cN_{d,e}^+$), and $\deg \tilde u> \deg \tilde v$
otherwise. In the former case the induction step goes as follows
(the proof in the latter case is similar). By lemma
\ref{lem-descr-normalis-bis}(a) we can write $\varphi_s$ as in
(\ref{eq-Jonq-plus}), where $f\in A_{d,e}$ and $\deg f\ge 2$.
Letting
$$\varphi:(x,y)\longmapsto (u,v)=(\alpha\tilde u+f(\tilde
v),\beta\tilde v)\,$$ from (\ref{eq-induct}) we obtain $$\deg
u=\deg f(\tilde v)=\deg f\cdot \deg\tilde v=\deg\varphi_s\cdot\deg
\psi=\prod_{i=1}^s\deg\varphi_i\ge 2\deg\tilde v>\deg\tilde v=\deg
v\,.$$ This ends the proof. \eproof

Now we can deduce the first part of proposition \ref{prop-amalga}.

\blem\la{lem-ge-1} If $e>1$ then \be\la{eq-isomm}
\langle\cN_{d,e}^+,\,\cN_{d,e}^-\rangle
\simeq\cN_{d,e}^+*_{\TT}\cN_{d,e}^-\,\,.\ee In particular, if
$e^2\not\equiv 1\mod d$ then \be\la{eq-isomm-bis}
\cN_{d,e}\simeq\cN_{d,e}^+*_{\TT}\cN_{d,e}^-\,.\ee \elem

\bproof The second assertion follows from the first by virtue of
lemma \ref{lem-generat}. To show the first one we recall that for
$e>1$
$$\cN_{d,e}^+\cap\cN_{d,e}^-=\TT=N^+_{d,e}=N^-_{d,e}\,,$$ see
(\ref{nde}) and (\ref{tor-int}). By a standard procedure (see
\cite[Ch. IX, \S 35, (6)]{Ku}) any element $\Phi\in
\langle\cN_{d,e}^+,\,\cN_{d,e}^-\rangle$ can be written in the
form
$$\Phi=t\varphi=t\varphi_s\cdot\ldots\cdot\varphi_1\,,$$ where
$t\in\TT$, $s\ge 0$, and for $s>0$ the factors
$\varphi_i\not\in\TT$ are as in (\ref{eq-alt-pro}). If  $s>0$ then
by (\ref{eq-deg-prod}) we obtain
$$\deg\Phi=\deg\varphi=\prod_{i=1}^{s}\deg\varphi_i>1\,.$$ Hence
$\Phi=\id$ if and only if $s=0$ and $t=\id$. Thus there is no
non-trivial relation in the group
$\langle\cN_{d,e}^+,\,\cN_{d,e}^-\rangle$ between elements of the
generating subgroups $\cN_{d,e}^\pm$ and so (\ref{eq-isomm}) holds
(cf.\ \cite[\S 13]{Wr3}).\eproof

To finish the proof of proposition \ref{prop-amalga} we need the
following auxiliary result from the combinatorial group theory.

\bthm\label{thm-Hanna} {\rm (Hanna Heumann \cite[Corollary
8.11]{Ne})} In the amalgamated free product $G=A*_C B$ with the
unified subgroup $C=A\cap B$, consider two subgroups $\tilde
A\subseteq A$ and $\tilde B\subseteq B$, and let $\tilde G=\langle
\tilde A,\tilde B\rangle$. Assume that $\tilde A\cap C=\tilde
C=\tilde B\cap C$. Then $\tilde G=\tilde A*_{\tilde C}\tilde
B$.\ethm

The next lemma proves the second part of proposition
\ref{prop-amalga}.

\blem\la{lem-ge-2} For $e^2\equiv 1\mod d$ we have
\be\la{eq-star-prod-bis-bis} \cN_{d,e}\simeq\cN_{d,e}^+
*_{N_{d,e}^+} N_{d,e}\,.\ee In particular, for $e=1$
\be\la{eq-star-prod} \cN_{d,1}\simeq\cN_{d,1}^+ *_{B^+}
\GL(2,\CC),\quad \text{where}\quad B^+=\cN_{d,e}^+\cap
\GL(2,\CC)\,,\ee while if $e>1$ and $e^2\equiv 1\mod d$ then
\be\la{eq-star-prod-bis} \cN_{d,e}\simeq\cN_{d,e}^+ *_{\TT}
N(\TT),\quad \text{where}\quad \TT=\cN_{d,e}^+\cap N(\TT)\,. \ee
\elem

\bproof We assume in the sequel that $e^2\equiv 1\mod d$. Let us
note first that (\ref{eq-star-prod}) and (\ref{eq-star-prod-bis})
are formal consequences of (\ref{eq-star-prod-bis-bis}) since
\[ N_{d,e}=\begin{cases} \GL(2,\CC)
&\text{ if $e=1$}\\
N(\TT)&\text{ if $e>1$}\end{cases}\qquad\text{and}\qquad
N_{d,e}^+=\begin{cases}
B^+ &\text{ if $e=1$}\\
\TT&\text{ if $e>1\,$.}\end{cases}\] Let us show
(\ref{eq-star-prod-bis-bis}). It follows from our definitions and
lemma \ref{lem-stabil-axes} that \be\label{ints}
{\cN}^+_{d,e}\cap{\Aff}^+(\AA^2)=N_{d,e}\cap{\Aff}^+(\AA^2)=N^+_{d,e}\,.\ee
Letting $$ A={\JONQ}^+(\AA^2),\quad
B=\Aff(\AA^2),\quad\text{and}\quad C=A\cap B={\Aff}^+(\AA^2)\,$$
by the Jung-van der Kulk theorem \ref{thm-JvdK} we obtain
$G=\Aut(\AA^2)=A*_C B$. Letting further
$$\tilde A={\cN}^+_{d,e},\quad \tilde B=N_{d,e},\quad\text{and}\quad \tilde
C=N^+_{d,e}=\tilde A\cap\tilde B\,$$  we see by (\ref{ints}) that
$\tilde A\cap C=\tilde B\cap C=\tilde C$. Now
(\ref{eq-star-prod-bis-bis}) follows by applying Hanna Neumann's
theorem \ref{thm-Hanna}. \eproof

\subsection{Algebraic groups actions on affine toric surfaces}
We can deduce now the following analog of the Kambayashi-Wright
theorem \ref{thm-KW} for affine toric surfaces.

\bthm\la{thm-KW-toric} Let $G\subseteq\Aut(X_{d,e})$ be an
algebraic group acting on an affine toric surface $X_{d,e}$, where
as before $1\le e< d$ and $\gcd(d,e)=1$. \bnum\item[(a)] If
$e^2\not\equiv 1\mod d$ then $G$ is conjugate
 in the group $\Aut(X_{d,e})$ to a subgroup of one of the groups
 $\cN_{d,e}^+/G_{d,e}$ and $\cN_{d,e}^-/G_{d,e}$.
 \item[(b)]
 If $e^2\equiv 1\mod d$ and $e>1$ then $G$ is conjugate
 to a subgroup of  one of the groups $\cN_{d,e}^+/G_{d,e}$
 and $N(\TT)/G_{d,e}$.\item[(c)] If $e=1$ then $G$
 is conjugate to a subgroup of  one of the groups
 $\cN_{d,e}^+/G_{d,e}$ and $\GL(2,\CC)/G_{d,e}$.\enum
\ethm

\bproof Consider the canonical surjection
$$\pi_*:\cN_{d,e}\to\cN_{d,e}/G_{d,e}\simeq \Aut(X_{d,e})\,$$
(see lemma \ref{lem-nornalis}). The algebraic group $\tilde
G=\pi_*^{-1}(G)$ has bounded degree. Under the assumption of (a)
we can conclude that the length $s={\rm length}(\varphi)$ in
(\ref{eq-alt-pro}) is uniformly bounded for all $\varphi\in\tilde
G\setminus\TT$. The same holds for $\pi_*(\varphi)\in G$ with
respect to the free amalgamated product structure
(\ref{eq-amalg-1}) as in theorem \ref{thm-JvdK-toric-surf}. Now
(a) follows by Serre's theorem \cite{Se} (cf.\ subsection
\ref{ssec-Aut}).  Due to
(\ref{eq-star-prod-bis-bis})--(\ref{eq-star-prod-bis}) a similar
argument applies also in the remaining cases (b) and (c). \eproof

The following corollary is immediate.

\bcor\la{cor-unipo} Any connected, unipotent subgroup
$U\subseteq\Aut(X_{d,e})$ is abelian.
 \ecor

For a reductive group acting on $X_{d,e}$,  the following analog
of proposition \ref{prop-quasitorus} holds.

\bthm\la{thm-rdgrac} Let $G\subseteq\Aut(X_{d,e})$ be a reductive
algebraic group. \bnum\item[(a)] If $e^2\not\equiv 1\mod d$ then
$G$ is conjugate
 in the group $\Aut(X_{d,e})$ to a subgroup of the torus $\TT/G_{d,e}$.
\item[(b)] If $e^2\equiv 1\mod d$ and $e>1$ then $G$ is conjugate
to a subgroup of  the quotient $N(\TT)/G_{d,e}$.
\item[(c)] If $e=1$ then $G$ is conjugate to a subgroup of
the quotient $\GL(2,\CC)/G_{d,e}$.\enum \ethm

\bproof Clearly the group $\tilde
G=\pi_*^{-1}(G)\subseteq\cN_{d,e}$ is reductive. By theorem
\ref{thm-KW-toric} we may assume that $\tilde G$ is a subgroup of
one of the corresponding factors. It suffices to restrict to the
case where $\tilde G\subseteq\cN_{d,e}^\pm$, since in the other
case the assertions are evidently true. In the former case by
proposition \ref{prop-quasitorus}  the group $\tilde G$ is abelian
and conjugate to a subgroup of the torus $\TT$ via an element
$\mu\in U^\pm$. We claim that such an element $\mu$ can be chosen
within the subgroup $U^\pm\cap\cN_{d,e}^\pm$. Let us show this
assuming that $\tilde G\subseteq\cN_{d,e}^+$, the other case being
similar. Indeed, consider an element $\varphi^+\in \tilde
G\subseteq \cN_{d,e}^+$ as in (\ref{eq-Jonq-plus}). By lemma
\ref{lem-descr-normalis-bis} we have  $f=\sum_{m\ge 0} a_my^m\in
A_{d,e}$, that is, $a_m=0$ $\forall m\not\equiv e\mod d$. Hence in
(\ref{eq-syst}) we can also choose $g=\sum_{m\ge 0}
b_my^m\in\kk[y]$ so that $b_m=0$ $\forall m\not\equiv e\mod d$ and
so $g\in  A_{d,e}$ too. Again by  lemma
\ref{lem-descr-normalis-bis} the latter ensures that $\mu\in
U^+\cap\cN_{d,e}^+$. Since $\tilde G$ is abelian, by virtue of
lemma \ref{lem-commuting} such an element $\mu$ can be found
simultaneously for all $\varphi^+\in \tilde G$. Since $\mu$
normalizes the cyclic group $G_{d,e}$ it descends to an
automorphism $\bar\varphi\in\Aut(X_{d,e})$ that conjugates $G$ to
a subgroup of the torus $\TT/G_{d,e}$. The proof is completed.
\eproof

The following simple example clarifies case (c) above.

\bexa\la{rem-ver-cone} The surface $X_{d,1}$ is a Veronese cone
i.e., the affine cone over a rational normal curve $\Gamma_d$ of
degree $d$ in $\PP^d$. The group $\GL(2,\CC)$ acts naturally on
$X_{d,1}$ (inducing the action of ${\rm PGL}(2,\CC)$ on
$\Gamma_d$) via the standard irreducible representation on the
space of binary forms of degree $d$. \eexa


\section{Acyclic curves and automorphism groups of non-toric quotient surfaces}
\label{sec-non-toric}

In this section we classify acyclic curves on the quotient
$X=\AA^2/\, G$ of the affine plane by a nonabelian finite linear
group $G$ and describe the automorphism groups of such surfaces.

We assume in the sequel that the finite group $G\subseteq \GL(n,
\CC)$ is small that is, does not contain any pseudoreflection.
Recall that a {\it pseudoreflection} on $\AA^n$ is a non-identical
linear transformation of finite order fixing pointwise a
hyperplane.  By the Chevalley-Shephard-Todd Theorem, the quotient
space $\AA^n/\, G$ of a finite linear group $G\subseteq\GL(n)$ is
isomorphic to $\AA^n$  if and only if $G$ is generated by
pseudoreflections. Assuming that this is not the case and
considering the normal subgroup $H\lhd\, G$ generated by all
pseudoreflections, we can decompose the quotient morphism $\pi
\colon \AA^n \to \AA^n/\, G$ into a two-step factorization
$$\AA^n\to\AA^n/H\simeq\AA^n\to X=\AA^n/G\,.$$ For $n=2$
we obtain in this way a presentation $X\simeq \AA^2/(G/H)$, where
the small linear group $G/H$ acts on $\AA^2$ freely off the
origin. Thus we may assume, without loss of generality, that $G$
is small. Under this assumption $X=\AA^2/\, G$ is a normal affine
surface with a unique singular point $Q=\pi(\bar 0)$, and the
quotient morphism $\pi : \AA^2\to X$ is unramified outside this
point.

It  is well known  that any finite subgroup of the group
$\SL(2,\CC)$ either is cyclic, or is a binary dihedral
(respectively, binary tetrahedral, octahedral, icosahedral) group,
see e.g. \cite{CS} or \cite{Ko}. The finite subgroups of
$\GL(2,\CC)$ are cyclic extensions of these
 groups; we refer  to
\cite{NvdPT} for their description.

For $n=2$ every small abelian group $G$ is conjugate to a cyclic
group $G_{d,\, e}$ and so $X=\AA^2/G$ is a toric surface. If $G$
is  small but nonabelian then $X$ is non-toric. In the next
theorem we examine this alternative possibility.

\bthm \label{prop6} Let $X=\AA^2/\, G$, where $G \subseteq
\GL_2(\CC)$ is a  nonabelian, finite, small
 group. Then any irreducible acyclic curve $\Oc$
in $X$ is the image $\pi(L)$ of  an affine line $L\subseteq \AA^2$
passing through the origin. \ethm

\bproof Let $\Oc$ be an irreducible acyclic curve on $X$. Assume
first that $G$ does not pass through the singular point
$Q=\pi(\bar 0)\in X$. Arguing as in the proof of
theorem~\ref{prop3} we get
$$
\pi^{*}(\Oc) \, = \, C^1 +\ldots+ C^d\,,
$$
where $d=|G|$, the irreducible components $C^i$ are disjoint,
$\pi:C^i\to \Oc\simeq\AA^1$ is an isomorphism for all $i$, and the
group $G$ acts transitively on the set of these components. Since
$d>1$, in suitable new coordinates $(x,y)$ on $\AA^2$ the curves
$C^i$ are given by equations $y=\kappa_i$ with
$\kappa_i\in\CC^\times$. Any $g\in G$ sends $y$ to $\lambda_g y$,
where $g\mapsto\lambda_g$ defines a character $\psi \colon G \to
\G_m$. The $G$-action on the set of components being simply
transitive, $\psi$ is injective and so $G$ must be abelian,
contrary to our assumption. This contradiction shows that
$Q\in\Oc$.

Consider the reduced acyclic plane curve
$C=\pi^{*}(\Oc)\subseteq\AA^2$ passing through the origin. If $C$
were smooth at the origin then the finite linear group $G$
preserving $C$ would preserve also the tangent line to $C$ at the
origin. So by Maschke's Theorem $G$ should be abelian, a
contradiction. Thus $\bar 0\in C$ is a singular point.

By corollary \ref{final-cor} the stabilizer $\Stab(C)\supseteq G$
is abelian unless $C$ is equivalent to a union of $r\ge 2$ affine
lines through the origin. Since $G$ is assumed to be nonabelian,
the latter condition holds indeed.

We claim that $C$ is in fact a union of lines, which yields the
assertion. Indeed, consider the  irreducible decomposition
$C=C^1+\ldots+C^r$, where $r\ge 2$. Notice that every component
$C^i$ of $C$ has a unique place at infinity. Since for every
$j\neq i$, $C^j=g(C^i)$ for some $g\in G$, where $g$ is linear, we
have $\deg(C^j)=\deg(C^i)=\delta$. Assume to the contrary that
$\delta>1$. Since $G$ is nonabelian, by Maschke's Theorem the
tautological representation $G\hookrightarrow\GL(2,\CC)$ is
irreducible. Hence the induced $G$-action on $\PP^1=\PP(\AA^2)$
has no fixed point. It follows that for some pair of indices
$i\neq j$, the points at infinity of the projective curves
$\overline{C^i}$ and $\overline{C^j}=\overline{g(C^i)}$ are
distinct. By Bezout Theorem these curves meet in $\delta^2>1$
points in $\PP^2$. Since all these points are situated in the
affine part $\AA^2$ and the intersection index of $C^i$ and $C^j$
at the origin equals $1$, these curves must have extra
intersection points in the affine part. This contradicts the fact
that $C$ is simply connected. This contradiction ends the proof.
\eproof

The following corollary  is immediate. Part (a) is an analog of
theorem~\ref{tmain}.

\bcor\label{cor-non-abel-orbits} Under the assumptions of
theorem~\ref{prop6} the following hold. \bnum\item[(a)] The only
irreducible acyclic curves in $X$ are the closures of
one-dimensional orbits of the $\G_m$-action on $X$ by
homotheties.\item[(b)]  The only simply connected curves in $X$
are the images of finite unions of affine lines through the origin
in $\AA^2$.\enum \ecor

Using theorem~\ref{prop6} we can deduce the following description
of the automorphism group of $X$, and as well an information on
the equivalence classes of irreducible acyclic curves in $X$.

\bthm\label{thm-non-abel-aut} Let as before $X=\AA^2/G$, where $G$
is a nonabelian small finite subgroup of $\GL(2,\CC)$. Then the
following hold. \bnum
\item[(a)]
$\Aut(X)\simeq N(G)/G\,,$ where the normalizer $N(G)$ of $G$ in
$\GL(2,\CC)$ is a finite extension of the one-torus
$\TT_{1,1}=Z(\GL(2,\CC))$. Consequently, $\Aut(X)$ is a finite
extension of the one-torus $\TT_{1,1}/Z(G)$.
\item[(b)] Let $\cS$ denote the set of all irreducible acyclic curves
$\Oc=\pi(L)$ on $X$. Then the $\Aut(X)$-action on $\cS$ has finite
orbits, and the orbit space $\cS/\Aut(X)$  is a rational
curve.\enum \ethm

\bproof (a) We claim that every automorphism $\alpha$ of the
quotient surface $X=\AA^2/G$ comes from an automorphism of $\AA^2$
normalizing the group $G$. Indeed, since $G$ is small, the
restriction $\pi|_{\AA^2\setminus\{\bar 0\}}: \AA^2\setminus\{\bar
0\}\to X\setminus\{Q\}$ is an unramified Galois cover with the
Galois group $G$. Since the surface $\AA^2\setminus\{\bar 0\}$ is
simply connected, this is a universal cover. Furthermore,
$\alpha\in\Aut(X)$ fixes the singular point $Q$ and induces an
automorphism of the smooth locus $X\setminus\{Q\}$. By the
monodromy theorem, both compositions
$\alpha\circ\pi,\,\alpha^{-1}\circ\pi:\AA^2\setminus\{\bar 0\}\to
X\setminus\{Q\}$ can be lifted to holomorphic maps $\tilde \alpha,
\widetilde {\alpha^{-1}}:\AA^2\setminus\{\bar
0\}\to\AA^2\setminus\{\bar 0\}$ in such a way that $\tilde
\alpha\circ \widetilde {\alpha^{-1}}=\id_{\AA^2\setminus\{\bar
0\}}$. According to the Hartogs Principle, the holomorphic
automorphism $\tilde \alpha$ of $\AA^2\setminus\{\bar 0\}$ extends
to such an automorphism of the whole plane $\AA^2$. Let us show
that $\tilde \alpha$ is (bi)regular. The induced homomorphism
$\tilde\alpha^*:\CC[x,y]\to {\rm Hol}(\AA^2)$ into the algebra of
entire holomorphic functions in two variables sends the ring of
invariants $\CC[x,y]^G$ to the polynomial ring $\CC[x,y]$. Hence
the entire holomorphic functions $f=\tilde\alpha^*(x)$ and
$g=\tilde\alpha^*(y)$ are integral over $\CC[x,y]$. Consequently,
they are of polynomial growth and so are polynomials. Since
$\tilde \alpha$ covers $\alpha$ it belongs to the normalizer
$\On(G)$ of $G$ in the group $\Aut(\AA^2)$. This proves our claim.

By theorem \ref{prop6} the image $\alpha(\Oc)$ of an irreducible
acyclic curve $\Oc=\pi(L)$ on $X$ is again such a curve. Hence any
lift $\tilde\alpha\in\Aut(\AA^2)$ of $\alpha$ preserves the
collection of lines through the origin. Similarly as in the proof
of corollary \ref{cross-lin}, it is easily seen that an
automorphism with the latter property is linear i.e.,
$\tilde\alpha\in\GL(2,\CC)$. Hence $\tilde\alpha\in N(G)$. So
finally $\Aut(X)=N(G)/G$. The proof of the remaining assertions is
easy and can be left to the reader. \eproof

\brems\label{normalizer} 1. Let as before $\On(G)$ denote the
normalizer of $G$ in the full automorphism group $\Aut(\AA^2)$.
The quotient group $\On(G)/G$ embeds into $\Aut(X)$. By virtue of
theorem~\ref{thm-non-abel-aut} it follows that $\On(G)=N(G)$ for
any nonabelian small  subgroup $G\subseteq\GL(2,\CC)$.

2. In contrast, for a toric surface $X_{d,e}=\AA^2/G_{d,e}$ the
group $\Aut(X_{d,e})$ acts infinitely transitively on the smooth
locus $X_{\rm reg}$ i.e., $m$-transitively for all $m\ge 1$, see
\cite{AKZ}. In particular, the group $\Aut(X)$ is infinite
dimensional. The equivalence class of any curve $\Oc$ on $X$ under
the $\Aut(X)$-action is infinite dimensional, too. Indeed, given
an arbitrary set $K$ of $m$ distinct points on $X_{\rm reg}$, due
to the $m$-transitivity we can interpolate $K$ by an image of
$\Oc$ under a suitable automorphism (cf.\ corollary 4.18 in
\cite{AFKKZ}). \erems

\bexa Consider the quaternion group
$$
Q_8=\left\{\pm I_2, \pm
\begin{pmatrix}
i & 0 \\
0 & -i
\end{pmatrix},
\pm
\begin{pmatrix}
0 & 1 \\
-1 & 0
\end{pmatrix},
\pm
\begin{pmatrix}
0 & i \\
i & 0
\end{pmatrix}\right\}\subseteq\GL(2,\CC)\,.
$$ This is a non-splittable central extension of the Klein four-group
$V_4\simeq\ZZ/2\ZZ\times\ZZ/2\ZZ$ by the center $\{\pm I_2\}$ of
$Q_8$. The algebra of invariants $\CC[x,y]^{Q_8}$ is generated by
the homogeneous polynomials
$$
f_1=x^4+y^4, \quad f_2=x^2y^2, \quad\text{and}\quad f_3=x^5y-xy^5
$$
satisfying a relation of degree $12$. The map $X=\AA^2/G\to \AA^3$
defined by these polynomials identifies $X$ with a surface of
degree 12 in $\AA^3$ given by a quasihomogeneous polynomial
equation with
 weights $(2,2,3)$. By theorem \ref{prop6} every
irreducible acyclic curve $\Oc$ on $X$ is the image of an affine
line $L$ on $\AA^2$ passing through the origin. Such a curve $\Oc$
is smooth if and only if the polynomial~$f_3$ vanishes on $L$, if
and only if $L=\CC v$, where $v$ is one of the vectors \be
\label{6l} (1,0), \quad (0,1), \quad (1,1), \quad (1,-1), \quad
(1,i), \quad (1,-i)\,. \ee Thus $\Oc$ coincides with one of the
corresponding affine lines on $X$ passing through the singular
point $Q$. In particular, these three lines are the only affine
lines on $X$ in $\AA^3$, and the image of any embedding
$\AA^1\hookrightarrow X$ coincides with one of them. These lines
are equivalent under the action on $X$ of the automorphism group
$\Aut(X)=N(Q_8)/Q_8$. Indeed, the vectors (\ref{6l}) belong to the
same orbit of the normalizer $N(Q_8)\subseteq\GL(2,\CC)$. \eexa


%
\end{document}